\documentclass[a4paper]{article}
 \pdfoutput=1
 \UseRawInputEncoding
\usepackage[textwidth=14.5cm]{geometry}
\usepackage{blindtext}
\usepackage{amsmath,amsfonts,amssymb,amsthm}
\usepackage{latexsym}
\usepackage{enumerate}
\usepackage{graphics}
\usepackage{epsfig}
\usepackage{mathdots}
\usepackage{caption}
\usepackage{subcaption}
\usepackage{setspace}
\usepackage{multirow}
\usepackage{dcolumn}
\usepackage{url}
\usepackage{tikz,pgfplots}
\usetikzlibrary{arrows,backgrounds}
\usepackage{ragged2e}
\usepackage[numbers,sort&compress]{natbib}

\newtheorem{thm}{Theorem}

\newtheorem{defn}{Definition}  
\newtheorem{prop}{Proposition}  
\newtheorem{cor}{Corollary}  
\newtheorem{exmp}{Example}  
\newenvironment{solution}{\emph{Solution.}}  

\author{Yongwen Zhu\\
\small School of Mathematics and Information Science, Yantai University,\\Yantai City 264005, P.R. China
\\ \small Email: {zyw@ytu.edu.cn}
}

\title{On the Nine-Palace Arithmetic\\ ---A New Method of Mental Calculation
\footnote{Supported by National Natural Science Foundation of China (11771375)}
}

\date{}

\begin{document}

\begin{sloppypar}

\maketitle

\abstract{
Based on the nine-palace diagram, we establish the systematical geometric theory of arithmetic, which can realize the arithmetical addition, subtraction, multiplication, division and other operations thoroughly in the mind.
In this paper, we give a brief introduction to this theory, including the rotation invariance theorem, the vector addition and the lattice addition, the summation by the center of gravity, the law of multiplication ones and carries on the nine-palace diagram, the counting method, and so on.
This systematical arithmetic can be viewed as the mathematics based on the Chinese Luoshu Diagram so that it is not only of mathematical significance, but also of philosophical significance.
}


\textbf{MR(2010) Subject Classification} 11A99, 20N99

\section{Introduction}

On August 15, 20-year-old Neelakantha Bhanu Prakash of India became the first Asian and the first non-European in 23 years to win the world mental arithmetic gold medal, beating 29 competitors from 13 countries.
Journalist Manish Pandey described Prakash as the Usain Bolt of mathematics, which shows how fast he can do mental arithmetic.

Speed calculation or mental arithmetic, is the advanced part of arithmetic calculation ability, has an important role to the development of memory, thinking rigor and innovation ability. There have been many scholars in different fields to study this question from different viewpoints such as pedagogy, psychology, physiology and even physical science, see for example \cite{ClearmanJ,Gowers,MariaG,MassonN,ChenFeiyan-3,ChenFeiyan-4,yindushuxue,yindufeituo,yindususuan,ZhangHongjiang,lishijie,ZhangJingzhong,yuwenyongquan-1}.
\par
Relatively speaking, the most systematic theory of mental computation, is the speed calculation method of Shi Feng-shou \cite{shifengshou}, which key content is the 26 carry formulae used for multiplication, but to master and use these formulae skillfully, it is difficult without a long time of training. We improve this theory and establish the theory of shear products and small carry method, which is another set of more effective and systematic theory of speed calculation after Shi Feng-shou quick calculation method \cite{jiandaoji-zyw}.

This article will introduce the nine-palace arithmetic theory, which is an another set of fast calculation theory we create, and it has a very close relationship with the traditional Chinese culture, because this arithmetic theory is deeply inspired by the Hetu Diagram, Luoshu Diagram, Eight-Gua diagram and Nine-Palace diagram, etc., and is actually a mathematics based on the nine-palace diagram.

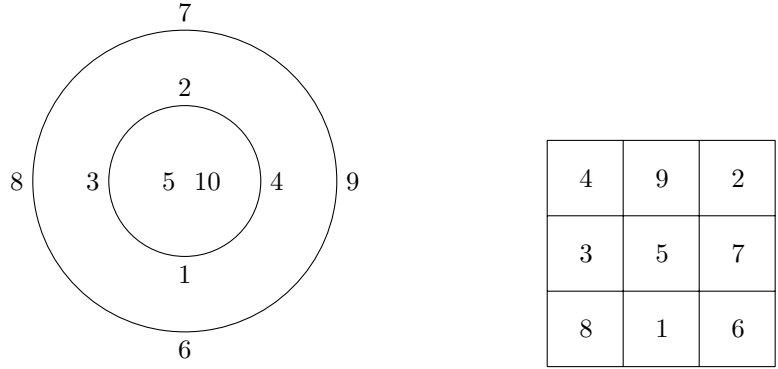
\begin{figure}[!htbp]
\centering
\begin{subfigure}[b]{.4\linewidth}
\centering
\begin{tikzpicture}
\draw (3,3) circle (1) (3,3) circle (2);
\coordinate [label=180:$5$](J) at (3,3);
\coordinate [label=0:$10$](J) at (3,3);
\coordinate [label=-90:$1$](J) at (3,2);
\coordinate [label=90:$2$](J) at (3,4);
\coordinate [label=180:$3$](J) at (2,3);
\coordinate [label=0:$4$](J) at (4,3);
\coordinate [label=-90:$6$](J) at (3,1);
\coordinate [label=90:$7$](J) at (3,5);
\coordinate [label=180:$8$](J) at (1,3);
\coordinate [label=0:$9$](J) at (5,3);
\end{tikzpicture}
\caption{Figures' Arrangement in Hetu Diagram}
\label{pic-HLshuzi-a}
\end{subfigure}
\quad
\centering
\begin{subfigure}[b]{.4\linewidth}
\centering
\begin{tikzpicture}
\draw[step=1] (0,0) grid (3,3);
\coordinate [label=0:$8$](A) at (0.3,0.5);
\coordinate [label=0:$1$](A) at (1.3,0.5);
\coordinate [label=0:$6$](A) at (2.3,0.5);
\coordinate [label=0:$3$](A) at (0.3,1.5);
\coordinate [label=0:$5$](A) at (1.3,1.5);
\coordinate [label=0:$7$](A) at (2.3,1.5);
\coordinate [label=0:$4$](A) at (0.3,2.5);
\coordinate [label=0:$9$](A) at (1.3,2.5);
\coordinate [label=0:$2$](A) at (2.3,2.5);
\end{tikzpicture}
\caption{Figures' Arrangement in the  Luoshu Diagram}
\label{pic-HLshuzi-b}
\end{subfigure}
\caption{Figures' Arrangement in Hetu and Luoshu diagrams}
\label{pic-HLshuzi}
\end{figure}

The Hetu Diagram and Luoshu Diagram is the origin of Chinese culture and an important part of ancient Chinese philosophy. It is an astronomical and geographical azimuth maps, a mathematical model map of the five elements in the universe and the development and change of everything.
The number arrangement of Hetu Diagram and Luoshu Diagram is particularly exquisite, which has its profound philosophical thought and also presents infinite beautiful mathematical rules.  See Fig.~\ref{pic-HLshuzi}.
Figures one and six in the Hetu Diagram are all in one group; Two, seven are with the same way; Three, eight are friends; Four, nine are friends too; Five and ten go the same way.
This means that each of the creating numbers is exactly 5 off the corresponding created numbers, and they are placed in the same direction, with 5 and 10 in the center. This is consistent with the principle of the ascension for 5 and the carry for 10 on the abacus. In this sense, the principle of the abacus is based on the Hetu Diagram.

If we divide the numbers $1,\ldots,9$ into three groups, that is, the group of small numbers $1,2,3$, that of middle numbers $4,5,6$, and that of large numbers $7,8,9$, then we see that each row and each column of the $3\times 3$ array in the Luoshu Diagram has exactly one small number, one middle number and one large number, respectively. If you calculate the $3\times 3$ determinant \cite{wangefang}, you will get exactly $360$, which corresponds to the approximate number of the all days of one year:

\begin{gather*}
\left|\begin{array}{cccc}
    4 & 9 & 2  \\ 3 & 5 & 7  \\8 & 1 & 6
\end{array}\right|   \\
=1\times 2\times 3+4\times 5\times 6+7\times 8\times 9-1\times 4\times 7-2\times 5\times 8-3\times 6\times 9=360.
\end{gather*}
Perhaps it is because the numbers are well matched in the order of large numbers, middle numbers and small numbers that the sum of the numbers in the $3\times 3$-array in each row, each column and each diagonal is equal to $15$. An $n\times n$ array is mathematically called an $n$-order \emph{magic square} if the sums of the numbers in each row, each column  and each diagonal are all equal, which is an important subject in combinatorics. Indeed, Luo Shu diagram is actually a magic square of the third order.
It is of great significance to further reveal and skillfully use the numerical rules in the Luoshu Diagram for the rapid calculation of arithmetic. As mentioned above, since the sum of the collinear three numbers is 15, the sum of the numbers on the intersecting lines is equal to 30 minus the common number on two lines; And since the sum of the two opposite numbers is 10 (complementary), the sum of the five numbers on two lines is equal  to 20 plus the complement of the common number on two lines. Now that 4,9,2 are collinear, 9,5,1 are collinear, and 9 is at the intersection of the two lines and its opposite point is 1,  we immediately get
\[4 + 9 + 2 + 5 + 1 = 21,\]
which shows a speed calculation.

The nine-palace arithmetic that we have created are essentially based on the Luoshu Diagram. However, in order to adapt to modern life habits, we arrange the numbers according to the natural orders of natural numbers in the nine-palace diagram, which is the way of most computer number keys and telephone number keys.
This arrangement is different from that of numbers in the Luoshu Diagram, but many important rules are the same for two arrangement ways when applied to arithmetic. For example, odd numbers and even numbers intersect, opposite points are complementary with each other, the numbers of the heaven run clockwise and the numbers of the earth run counterclockwise, the rotation invariance property, the formula of the scarecrow and so on are all the same.
A series of important discoveries of us such as the rotation invariance principle for addition and multiplication, the straw man formula, and so on, makes our nine-palace arithmetic theory completely distinguished from all other finger algorithms including the ancient Chinese calculation method of the swallowing gold in the sleeve.
Since it is an arithmetical theory based on the Luoshu Diagram, this theory has not only mathematical significance but also philosophical significance.

\section{Preliminaries: the nine-palace diagram and its extension}

We have solved the classical problem of the nine-palace rearrangement \cite{Zhu-8} by means of group theory.
The so-called \emph{basic nine-palace grid} is a table with three rows and three columns, each of which is a square, called a \emph{palace}, as shown in  Fig.~\ref{jiugongge-jiben-a}.
Now, we use the nine-palace grids for mental arithmetic.
The so-called nine-palace grid mental arithmetic method is to use nine-palace grid to complete arithmetic operations, and its basic idea is to use the position of the palace to represent the numbers. Fig.~\ref{jiugongge-jiben-b} is called the \emph{primitive nine-palace diagram}, in which each one of the nine-palace in turn represents numbers $1$, $2$, $3$, $4$, $5$, $6$, $7$, $8$ and $9$. Note that the core idea here is to use the palace to represent the number so that you don't need to actually fill the number in the palace as shown here. But if we want to do a nine-palace grid calculation in our head, then we first have to remember the corresponding relationship between numbers and positions, that is to say, we must memorize Fig.~\ref{jiugongge-jiben-b} proficiently.

\begin{figure}[!htbp]
\centering
\begin{subfigure}[b]{.4\linewidth}
\centering
\begin{tikzpicture}
\draw[step=1] (0,0) grid (3,3);
\end{tikzpicture}
\caption{The blank nine-palace diagram} \label{jiugongge-jiben-a}
\end{subfigure}
\quad
\begin{subfigure}[b]{.4\linewidth}
\centering
\begin{tikzpicture}
\draw[step=1] (0,0) grid (3,3);
\coordinate [label=0:$7$](A) at (0.3,0.5);
\coordinate [label=0:$8$](A) at (1.3,0.5);
\coordinate [label=0:$9$](A) at (2.3,0.5);
\coordinate [label=0:$4$](A) at (0.3,1.5);
\coordinate [label=0:$5$](A) at (1.3,1.5);
\coordinate [label=0:$6$](A) at (2.3,1.5);
\coordinate [label=0:$1$](A) at (0.3,2.5);
\coordinate [label=0:$2$](A) at (1.3,2.5);
\coordinate [label=0:$3$](A) at (2.3,2.5);
\end{tikzpicture}
\caption{The primitive nine-palace diagram}\label{jiugongge-jiben-b}
\end{subfigure}
\caption{The basic nine-palace grid}
\label{jiugongge-jiben}
\end{figure}

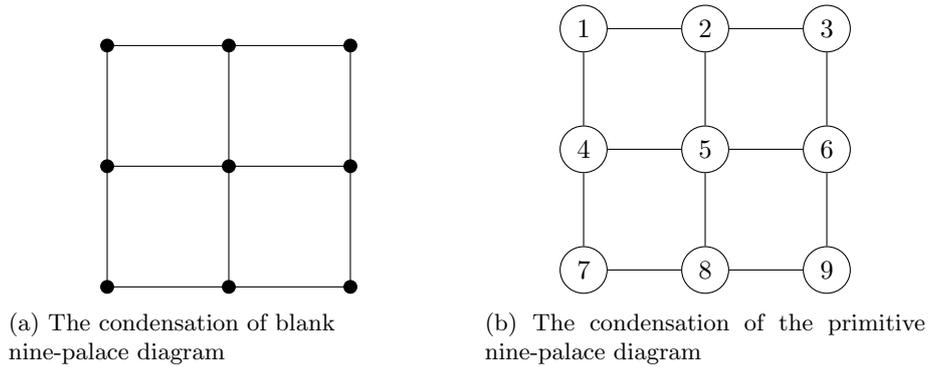
\begin{figure}[!htbp]
\centering
\begin{subfigure}[b]{.4\linewidth}
\centering
\begin{tikzpicture}[scale=0.8]
\draw[step=2] (0,0) grid (4,4);
\foreach \x in {0,2,4} \foreach \y in {0,2,4} \draw[fill] (\x,\y) circle (3pt);
\end{tikzpicture}
\caption{The condensation of blank \\nine-palace diagram} \label{benyuanjibenjiugongtu-a}
\end{subfigure}
\quad
\begin{subfigure}[b]{.4\linewidth}
\centering
\begin{tikzpicture}[scale=0.8]
\node(1) at (0,4) [circle,draw]{1};
\node(2) at (2,4) [circle,draw]{2};
\node(3) at (4,4) [circle,draw]{3};
\draw (1)--(2)--(3);
\node(4) at (0,2) [circle,draw]{4};
\node(5) at (2,2) [circle,draw]{5};
\node(6) at (4,2) [circle,draw]{6};
\draw (4)--(5)--(6);
\node(7) at (0,0) [circle,draw]{7};
\node(8) at (2,0) [circle,draw]{8};
\node(9) at (4,0) [circle,draw]{9};
\draw (7)--(8)--(9);
\draw (1)--(4)--(7);
\draw (2)--(5)--(8);
\draw (3)--(6)--(9);
\end{tikzpicture}
\caption{The condensation of the primitive nine-palace diagram} \label{benyuanjibenjiugongtu-b}
\end{subfigure}
\caption{The condensation of nine-palace diagram} \label{benyuanjibenjiugongtu}
\end{figure}

When we actually draw a nine-palace diagram, we usually draw each palace by condensing it into a point. The two sub-figures in Fig.~\ref{jiugongge-jiben} become two sub-figures in Fig.~\ref{benyuanjibenjiugongtu} after condensation. In this way, the condensed nine-palace grids have actually become the field grid,  but we still call them nine-palace grids or nine-palace diagrams.

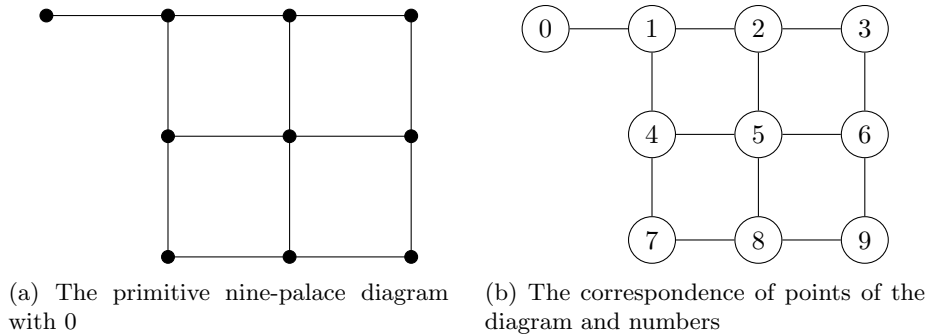
\begin{figure}[!htbp]
\centering
\begin{subfigure}[b]{.4\linewidth}
\centering
\begin{tikzpicture}[scale=0.8]
\draw[step=2] (0,0) grid (4,4);
\foreach \x in {0,2,4} \foreach \y in {0,2,4} \draw[fill] (\x,\y) circle (3pt);
\draw[fill] (-2,4) circle(3pt);
\draw (-2,4)--(0,4);
\end{tikzpicture}
\caption{The primitive nine-palace diagram with 0} \label{jitongtuwith0-a}
\end{subfigure}
\quad
\begin{subfigure}[b]{.4\linewidth}
\centering
\begin{tikzpicture}[scale=0.7]
\node(0) at (-2,4) [circle,draw]{0};
\node(1) at (0,4) [circle,draw]{1};
\node(2) at (2,4) [circle,draw]{2};
\node(3) at (4,4) [circle,draw]{3};
\draw (0)--(1)--(2)--(3);
\node(4) at (0,2) [circle,draw]{4};
\node(5) at (2,2) [circle,draw]{5};
\node(6) at (4,2) [circle,draw]{6};
\draw (4)--(5)--(6);
\node(7) at (0,0) [circle,draw]{7};
\node(8) at (2,0) [circle,draw]{8};
\node(9) at (4,0) [circle,draw]{9};
\draw (7)--(8)--(9);
\draw (1)--(4)--(7);
\draw (2)--(5)--(8);
\draw (3)--(6)--(9);
\end{tikzpicture}
             \caption{The correspondence of points of the diagram and numbers} \label{jitongtuwith0-b}
\end{subfigure}
\caption{The primitive nine-palace diagram with 0}
\label{jitongtuwith0}
\end{figure}

We note that the number $0$ is currently ignored. If we want the number $0$ to be represented in the diagram, we need to extend the graph. If the original nine-palace diagram is expanded into the diagram as in Fig.~\ref{jitongtuwith0}, it is called the \emph{primitive nine-palace diagram with zero}, or it still be called \emph{primitive nine-palace diagram} for short. That is to say, the primitive nine-palace diagram we mentioned later can be with or without zero.

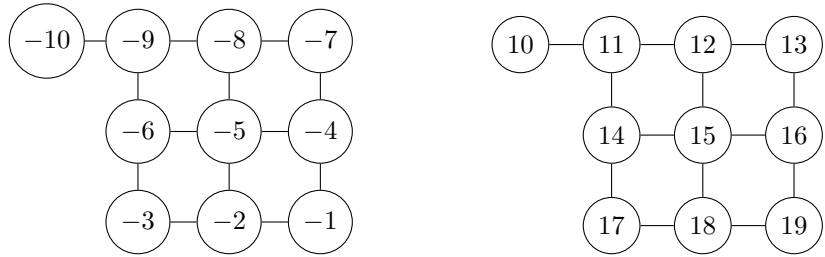
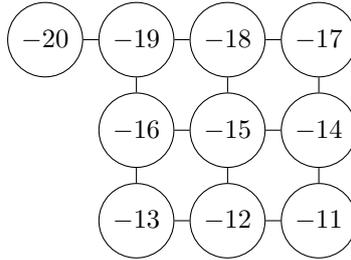
\begin{figure}[!htbp]
\centering
\begin{subfigure}[b]{.4\linewidth}
\centering
\begin{tikzpicture}[scale=0.6]
\node(0) at (-2,4) [circle,draw]{$-10$};
\node(1) at (0,4) [circle,draw]{$-9$};
\node(2) at (2,4) [circle,draw]{$-8$};
\node(3) at (4,4) [circle,draw]{$-7$};
\draw (0)--(1)--(2)--(3);
\node(4) at (0,2) [circle,draw]{$-6$};
\node(5) at (2,2) [circle,draw]{$-5$};
\node(6) at (4,2) [circle,draw]{$-4$};
\draw (4)--(5)--(6);
\node(7) at (0,0) [circle,draw]{$-3$};
\node(8) at (2,0) [circle,draw]{$-2$};
\node(9) at (4,0) [circle,draw]{$-1$};
\draw (7)--(8)--(9);
\draw (1)--(4)--(7);
\draw (2)--(5)--(8);
\draw (3)--(6)--(9);
\end{tikzpicture}
\caption{The nine-palace diagram of $-1$-family} \label{kzujiugongtu-a}
\end{subfigure}
\quad
\begin{subfigure}[b]{.4\linewidth}
\centering
\begin{tikzpicture}[scale=0.6]
\node(0) at (-2,4) [circle,draw]{10};
\node(1) at (0,4) [circle,draw]{11};
\node(2) at (2,4) [circle,draw]{12};
\node(3) at (4,4) [circle,draw]{13};
\draw (1)--(2)--(3);
\node(4) at (0,2) [circle,draw]{14};
\node(5) at (2,2) [circle,draw]{15};
\node(6) at (4,2) [circle,draw]{16};
\draw (4)--(5)--(6);
\node(7) at (0,0) [circle,draw]{17};
\node(8) at (2,0) [circle,draw]{18};
\node(9) at (4,0) [circle,draw]{19};
\draw (7)--(8)--(9);
\draw (0)--(1)--(4)--(7);
\draw (2)--(5)--(8);
\draw (3)--(6)--(9);
\end{tikzpicture}
\caption{The nine-palace diagram of $+1$-family} \label{kzujiugongtu-b}
\end{subfigure}
\quad
\begin{subfigure}[b]{.4\linewidth}
\centering
\begin{tikzpicture}[scale=0.6]
\node(0) at (-2,4) [circle,draw]{$-20$};
\node(1) at (0,4) [circle,draw]{$-19$};
\node(2) at (2,4) [circle,draw]{$-18$};
\node(3) at (4,4) [circle,draw]{$-17$};
\draw (1)--(2)--(3);
\node(4) at (0,2) [circle,draw]{$-16$};
\node(5) at (2,2) [circle,draw]{$-15$};
\node(6) at (4,2) [circle,draw]{$-14$};
\draw (4)--(5)--(6);
\node(7) at (0,0) [circle,draw]{$-13$};
\node(8) at (2,0) [circle,draw]{$-12$};
\node(9) at (4,0) [circle,draw]{$-11$};
\draw (7)--(8)--(9);
\draw (0)--(1)--(4)--(7);
\draw (2)--(5)--(8);
\draw (3)--(6)--(9);
\end{tikzpicture}
\caption{The nine-palace diagram of $-2$-family} \label{kzujiugongtu-c}
\end{subfigure}
\caption{The examples of the nine-palace diagrams of $k$-family} \label{kzujiugongtu}
\end{figure}

The significance of extending the nine-palace diagram is that we can go beyond the scope of the nine-palace in the process of using the original nine-palace diagram to calculate.
Let $k$ be an integer. If the number represented by each point in a primitive nine-palace diagram with zero is added $k$ multiples of $10$, then the resulting nine-palace diagram is called the \emph{$k$-family nine-palace diagram} (with the \emph{base point} $10k$), where $k$ is called the \emph{family number}.
Families with positive family numbers are collectively \emph{positive}, and families with negative family numbers are collectively \emph{negative}.
Fig.~\ref{kzujiugongtu} shows the nine-palace diagrams of $-1$-family, $+1$-family and $-2$-family respectively. The figure shows the natural connection between the two nine-palace diagram of the $-1$-family and that of the primitive family.

\begin{figure}[!htbp]
\centering
\begin{tikzpicture}[scale=0.7]
\node(01) at (2,10) [circle,draw]{$-9$};
\node(02) at (4,10) [circle,draw]{$-8$};
\node(03) at (6,10) [circle,draw]{$-7$};
\draw (01)--(02)--(03);
\node(04) at (2,8) [circle,draw]{$-6$};
\node(05) at (4,8) [circle,draw]{$-5$};
\node(06) at (6,8) [circle,draw]{$-4$};
\draw (04)--(05)--(06);
\node(07) at (2,6) [circle,draw]{$-3$};
\node(08) at (4,6) [circle,draw]{$-2$};
\node(09) at (6,6) [circle,draw]{$-1$};
\draw (07)--(08)--(09);
\draw (01)--(04)--(07);
\draw (02)--(05)--(08);
\draw (03)--(06)--(09);
\node(1) at (4,4) [circle,draw]{$1$};
\node(2) at (6,4) [circle,draw]{$2$};
\node(3) at (8,4) [circle,draw]{$3$};
\draw (1)--(2)--(3);
\node(4) at (4,2) [circle,draw]{$4$};
\node(5) at (6,2) [circle,draw]{$5$};
\node(6) at (8,2) [circle,draw]{$6$};
\draw (4)--(5)--(6);
\node(7) at (4,0) [circle,draw]{$7$};
\node(8) at (6,0) [circle,draw]{$8$};
\node(9) at (8,0) [circle,draw]{$9$};
\draw (7)--(8)--(9);
\draw (1)--(4)--(7);
\draw (2)--(5)--(8);
\draw (3)--(6)--(9);
\draw (08)--(1)(09)--(2);
\end{tikzpicture}
\caption{The conjunction of two nine-palace diagrams}
\label{+-benyuanduiqi}
\end{figure}
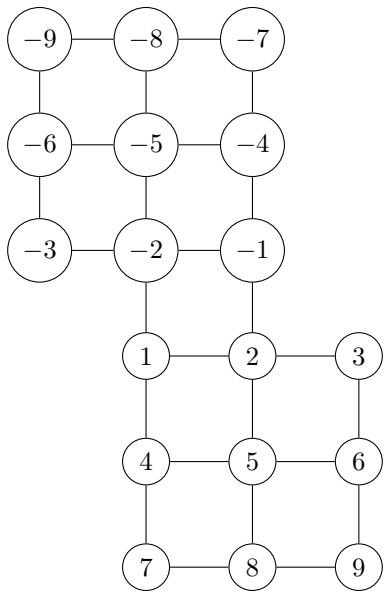

We read the \emph{position} in terms of numbers in the primitive nine-palace diagram with zero so that we can say position $0$, position $1$, position $2$, $\ldots$, position $9$. The position $10$ in a family is just the position 0 in the next family. Position $5$ is exactly the \emph{center} of the nine-palace diagram. Position $0$ is in position $1$ left away one space, position $10$ is in position $9$ right away one space.

Among the eleven positions on the nine-palace diagram, the eight points except the positions $0$, $5$, and $10$ just form a square, which we call the \emph{border} of the diagram. The border has four \emph{edges}, called the \emph{top, bottom, left, and right} edge according to their natural position respectively. The midpoint of the edges is also called the \emph{midpoint} of the nine-palace diagram. These are the \emph{up midpoint, bottom midpoint, left midpoint, right midpoint}, and the corresponding position is exactly $2,8,4,6$  respectively. The points in the corner of the frame are called the \emph{corner} of the nine-palace diagram. There are four corner points, i.e., the \emph{up left corner, up right corner, left bottom corner, bottom right corner}, which positions are $1,3,7,9$, respectively. Generally, we have

\begin{prop}
In the nine-palace diagram of $k$-family, the  position $n$ represents the number $k\times 10+n$, that is $(k, n) $. In particular, when $k < 0 $, the number is equal to $-(-k-1,10-n)$.
\end{prop}

As we have already seen above, with the extended stacking method, we can expand downwards the primitive nine-palace diagram to that of $-1$-family, and also we can expand the nine-palace diagram of $-1$-family upwards to the primitive one.
If you use the extended stacking method to continue to expand under the primitive nine-palace  diagram, then you could successively get the nine-palace  diagrams of 1--family, 2--family, 3--family, etc. Similarly,
above  the primitive nine-palace  diagram, you could successively get the nine-palace  diagrams of $-1$-family, $-2$-family, $-3$-family, and so on. It can be seen that the nine-palace diagrams of all families can be stacked together, with their base points on a slightly tilted vertical line.  The higher you go up, the smaller the numbers are; and as you go down, the numbers get more bigger.

\begin{figure}[!htbp] 
\centering
\begin{subfigure}[b]{.4\linewidth}
\centering
\begin{tikzpicture}[scale=0.6]
\draw[step=2] (0,0) grid (4,4);
\draw[fill] (0,4) circle (3pt)(4,4) circle (3pt)  (0,2) circle (3pt)  (2,2) circle (3pt)(4,0) circle (3pt);
\draw[line width=3pt,->,>=latex] (4,4)--(0,2)(0,2)--(2,2)--(4,0) ;
\draw[line width=3pt] (4,4) arc (0:180:2 and 1);
\draw[line width=3pt] (0,2) arc (0:360:0.5 and 0.5);
\coordinate [label=-90:$\quad$](A1) at (2,-0.5);
\end{tikzpicture}
\caption{$134\,
459$ \\in the primitive family}
\label{fig:duoweishubiaoshi-a}
\end{subfigure}
\quad
\begin{subfigure}[b]{.4\linewidth}
\centering
\begin{tikzpicture}[scale=0.6]
\draw[step=2] (0,0) grid (4,4);
\draw[fill] (2,0.2) circle (3pt)(2,0) circle (3pt)  (4,0) circle (3pt)(6,0) circle (3pt)  (2,2) circle (3pt)(4,2) circle (3pt);
\coordinate [label=135:$-1$](A) at (2,0);
\coordinate [label=225:$-1$](A1) at (2,0);
\coordinate [label=-90:$-1$](B) at (4,0);
\coordinate [label=-90:$-1$](C) at (6,0);
\coordinate [label=135:$-1$](D) at (2,2);
\coordinate [label=45:$-1$](E) at (4,2);
\draw[line width=3pt,->,>=latex](A)--(B)--(C)(C)--(D)(D)--(E);
\end{tikzpicture}
\caption{$-221\,054$ in the $-1$-family}
\label{fig:duoweishubiaoshi-b}
\end{subfigure}
\caption{Use a nine-palace diagram to represent long numbers}
\label{fig:duoweishubiaoshi}
\end{figure}
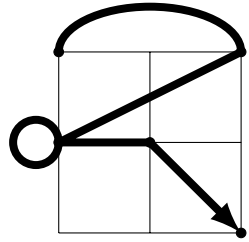
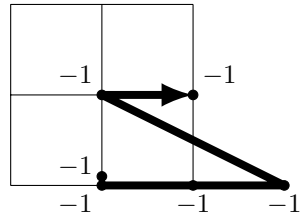

Speed calculation is inseparable from memory.
Matteo Ricci, an Italian missionary to China during the Wanli Period of the Ming Dynasty, wrote in Chinese a book titled The Memory Method of Western Countries, which was a famous work on the method of memorization \cite{lidouma}. The method of memorization described in this book is now commonly called a memory palace. We point out that the nine-palace diagram itself is a good kind of memory palaces. In order to remember numbers, we can use the nine-palace diagram. For example, $-1305$ is represented by points $-1$, $-3$, $0$ and $-5$ connected by some lines with an arrow indicating the order. Fig.~\ref{fig:duoweishubiaoshi} represents the positive integer $134\,459$ with the negative integer $-221\, 054$ clearly.

\section{Basic diagrams of addition and subtraction on the nine-palace diagram}

This section describes the basic diagrams for adding and subtracting on the nine-palace diagram, which represent how to add and subtract two numbers on the grid.\par

In the basic diagrams of addition, the arrow represents the method of adding a number, starting with the addend and ending with the result of adding the number. That is to say, starting from the point representing the addend in the grid, moving a distance of the length indicated by the arrow in the direction indicated by the arrow, the number corresponding to the point reached is the sum obtained. These diagrams are of course easy to understand, but the key is to remember the length and direction of the arrows in these basic diagrams, because this alone will speed up our calculations. These arrows essentially represent the vector \cite{qiuweisheng}.
The following is the illustration of the basic diagrams of  addends from 1 to 9.

\begin{figure}[!htbp] 
\centering
\begin{subfigure}[t]{.4\linewidth}
\begin{tikzpicture}[scale=0.6]
\draw[step=2] (0,0) grid (4,4);
\draw[fill]  (0,4) circle (3pt) (4,4) circle (3pt);
\draw[line width=3pt][->,>=latex] (0,4)--(4,4);
\end{tikzpicture}
\caption{The first picture of adding 2}\label{pic-plus2-a}
\end{subfigure}
\quad
\begin{subfigure}[t]{.4\linewidth}
\begin{tikzpicture}[scale=0.6]
\draw[step=2] (0,0) grid (4,4);
\draw[fill]  (2,4) circle (3pt) (0,2) circle (3pt);
\draw[line width=3pt] (2,4) [->,>=latex]--(0,2);
\end{tikzpicture}
\caption{The second picture of adding 2}\label{pic-plus2-b}
\end{subfigure}
\caption{The pictures of adding 2}
\label{pic-plus2}
\end{figure}
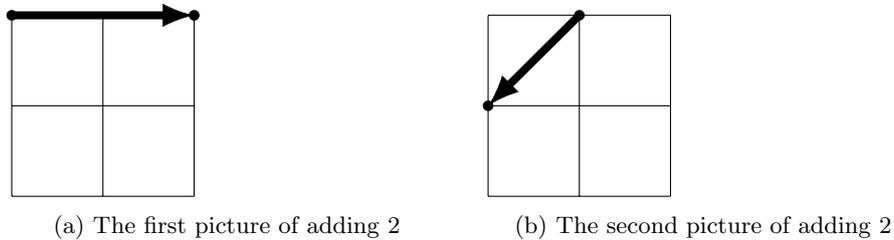

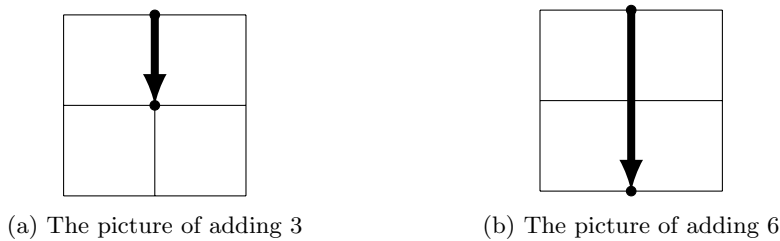
\begin{figure}[!htbp] 
\centering
\begin{subfigure}[t]{.4\linewidth}
\centering
\begin{tikzpicture}[scale=0.6]
\draw[step=2] (0,0) grid (4,4);
\draw[fill]  (2,4) circle (3pt) (2,2) circle (3pt);
\draw[line width=3pt] (2,4) [->,>=latex]--(2,2);
\end{tikzpicture}
\caption{The picture of adding 3}
\label{pic-plus3}
\end{subfigure}
\quad
\begin{subfigure}[t]{.4\linewidth}
\centering
\begin{tikzpicture}[scale=0.6]
\draw[step=2] (0,0) grid (4,4);
\draw[fill]  (2,4) circle (3pt) (2,0) circle (3pt);
\draw[line width=3pt] (2,4) [->,>=latex]--(2,0);
\end{tikzpicture}
\caption{The picture of adding 6}\label{pic-plus6}
\end{subfigure}
\caption{The pictures of adding 3 and adding 6}
\label{pic-plus-3-6}
\end{figure}

The way to add or subtract $k$ from a number is to move $k$ steps forwards or backwards in the grid.
As shown in Fig.~\ref{pic-plus1} the addition of a number plus $1$ is exhibited. Indeed, Fig.~\ref{pic-plus1-a}  represents  $7 + 1 = 8$ (in the primitive nine-palace diagram) or $8-3 + 1 = 2 $ (in the diagram of $-1$-family); and
 Fig.~\ref{pic-plus1-b} represents $3 + 1 = 4 $ (primitive family) or $7 + 1 = - $6 ($-1$-family).
For more examples, Figs.~\ref{pic-plus2} and~\ref{pic-plus-3-6} and so on represent the basic diagrams of addends $2$, $3$ and $6$ respectively. We're not going to draw the basic picture of adding the other numbers.

\begin{figure}[!htbp]  
\centering
\begin{subfigure}[t]{.4\linewidth}
\begin{tikzpicture}[scale=0.6]
\draw[step=2] (0,0) grid (4,4);
\draw[fill]  (0,0) circle (3pt) (2,0) circle (3pt);
\draw[line width=3pt][->,>=latex] (0,0)--(2,0);
\end{tikzpicture}
\caption{The first picture of adding 1}\label{pic-plus1-a}
\end{subfigure}
\quad
\begin{subfigure}[t]{.4\linewidth}
\begin{tikzpicture}[scale=0.6]
\draw[step=2] (0,0) grid (4,4);
\draw[fill]  (4,4) circle (3pt) (0,2) circle (3pt);
\draw[line width=3pt] (4,4) [->,>=latex]--(0,2);
\end{tikzpicture}
\caption{The second picture of adding 1}\label{pic-plus1-b}
\end{subfigure}
\caption{The pictures of adding 1}
\label{pic-plus1}
\end{figure}
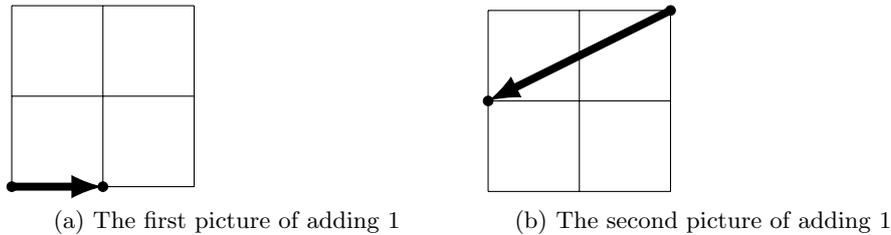

The basic diagrams of subtraction on the nine-palace diagram is the reverse sign of the arrow in the basic diagrams of addition.
Keeping in mind  the basic diagrams of addition and abstraction of all ten numbers is the basic skills of learning nine-palace arithmetic.

We point out that the addition between basic graphic arrows (called the vector addition in vector algebra \cite{wangefang,qiuweisheng}) satisfies a law called the \emph{triangle rule}. In accordance with the law,  the corresponding arrow of a number plus $m$ and then plus $n $  equals  the arrow of the number plus $ m + n$, this is because  $(a + m) + n = a + (m + n) $. Similarly, the corresponding arrow of  a number minus $m $ and then minus  $n $  equals the arrow of the number minus $m + n$, this is because  $(a - m) - n = a - (m + n) $. Adding and subtracting together satisfy the rule too. We write this rule as a theorem:

\begin{thm}\label{sanjiaoxingfaze}
Addition (subtraction) in the form of arrows on the nine-palace diagram satisfies the triangle rule.
\end{thm}

Since the numbers can be represented in the form of arrows on the grid, that is, as vectors, the triangle method of number operation on the graph is actually the triangle rule satisfied by the vectors' addition operation.

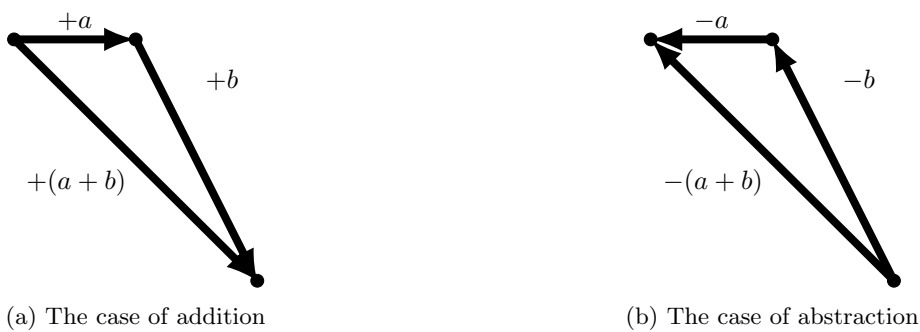
\begin{figure}[!htbp] 
\centering
\begin{subfigure}[t]{.4\linewidth}\centering
\begin{tikzpicture}[scale=0.8]
\draw[fill]  (0,4) circle (3pt) (4,0) circle (3pt) (2,4) circle (3pt);
\draw[line width=3pt][->,>=latex] (0,4)--(2,4);
\draw[line width=3pt][->,>=latex] (2,4)--(4,0);
\draw[line width=3pt][->,>=latex] (0,4)--(4,0);
\coordinate [label=90:$+a$](A) at (1,4);
\coordinate [label=45:$+b$](B) at (3,3);
\coordinate [label=225:$+(a+b)$](C) at (2,2);
\end{tikzpicture}
\caption{The case of addition} \label{triangle-law-a}
\end{subfigure}
\hskip70pt
\centering
\begin{subfigure}[t]{.4\linewidth}\centering
\begin{tikzpicture}[scale=0.8]
\draw[fill]  (0,4) circle (3pt) (4,0) circle (3pt) (2,4) circle (3pt);
\draw[line width=3pt][<-,>=latex] (0,4)--(2,4);
\draw[line width=3pt][<-,>=latex] (2,4)--(4,0);
\draw[line width=3pt][<-,>=latex] (0,4)--(4,0);
\coordinate [label=90:$-a$](A) at (1,4);
\coordinate [label=45:$-b$](B) at (3,3);
\coordinate [label=225:$-(a+b)$](C) at (2,2);
\end{tikzpicture}
\caption{The case of abstraction} \label{triangle-law-b}
\end{subfigure}
\caption{The triangle law} \label{triangle-law}
\end{figure}

The triangle rule is represented by a graph that happens to be a triangle, as shown in  Fig.~\ref{triangle-law}. For example, because $1+1=2$, the arrow for adding 2 is going in exactly the same direction as the arrow for adding 1, but the arrow for adding 2 is twice as long as that for 1. Similarly, the arrow of the addition plus 4 is in the same direction as the arrow of the addition plus 2, but the former is twice as long as the latter; The arrow of the addition plus 8 goes in exactly the same direction as that plus 4, but the former is twice as long as the latter.

\begin{figure}[!htbp] 
\centering
\begin{subfigure}[t]{.4\linewidth}\centering
\begin{tikzpicture}[scale=0.5]
\draw[step=2] (0,0) grid (4,4);
\draw[fill]  (2,4) circle (3pt) (2,2) circle (3pt)(0,2) circle (3pt)(-2,2) circle (3pt)(4,2) circle (3pt)(6,2) circle (3pt);
\draw[line width=3pt] (2,4) [->,>=latex]--(0,2);
\draw[line width=3pt] (2,4) [->,>=latex]--(-2,2);
\draw[line width=3pt] (2,4) [->,>=latex]--(2,2);
\draw[line width=3pt] (2,4) [->,>=latex]--(6,2);
\draw[line width=3pt] (2,4) [->,>=latex]--(4,2);
\coordinate [label=270:$+1$](A) at (-2,2);
\coordinate [label=270:$+2$](A) at (0,2);
\coordinate [label=270:$+3$](A) at (2,2);
\coordinate [label=270:$+4$](A) at (4,2);
\coordinate [label=270:$+5$](A) at (6,2);
\end{tikzpicture}
\caption{Symmetric with respect to $+3$}\label{duichen+3}
\end{subfigure}
\quad
\begin{subfigure}[t]{.4\linewidth}\centering
\begin{tikzpicture}[scale=0.5]
\draw[step=2] (0,0) grid (4,4);
\draw[fill]  (2,4) circle (3pt) (2,0) circle (3pt)(0,0) circle (3pt)(-2,0) circle (3pt)(4,0) circle (3pt)(6,0) circle (3pt);
\draw[line width=3pt] (2,4) [->,>=latex]--(0,0);
\draw[line width=3pt] (2,4) [->,>=latex]--(-2,0);
\draw[line width=3pt] (2,4) [->,>=latex]--(2,0);
\draw[line width=3pt] (2,4) [->,>=latex]--(6,0);
\draw[line width=3pt] (2,4) [->,>=latex]--(4,0);
\coordinate [label=270:$+4$](A) at (-2,0);
\coordinate [label=270:$+5$](A) at (0,0);
\coordinate [label=270:$+6$](A) at (2,0);
\coordinate [label=270:$+7$](A) at (4,0);
\coordinate [label=270:$+8$](A) at (6,0);
\end{tikzpicture}
\caption{Symmetric with respect to $+6$}\label{duichen+6}
\end{subfigure}
\caption{Symmetry relationship of diagrams for addition}\label{duichen+3+6}
\end{figure}
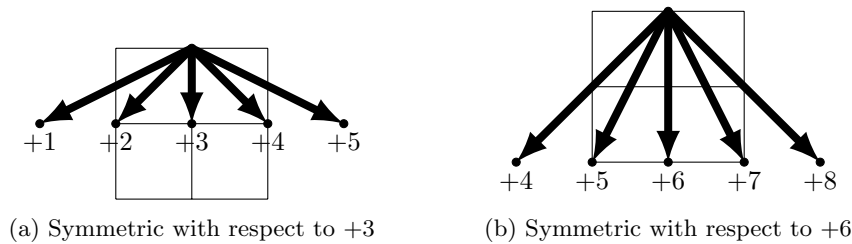

In addition, we notice some symmetric relationships between the basic figures. For example, the arrow of adding 3 is a vertical small vertical, and arrow of adding 2 and  the arrow of adding 4 are small prime and small stem, respectively. Two arrows of adding 2 and4 are symmetric with respect to that of adding 3, because $3-1=2$ and $3+1=4$. And similarly, the  arrows of adding 1 and 5 is also symmetric with respect to the arrow of adding 3. These symmetries are expressed by the graph ~\ref{duichen+3}.
The basic diagram for adding 6 is a long vertical. The arrows of adding 5 and 7, as well as  those of adding 4 and 8, are also symmetric with respect to the arrow of adding 6, see  Fig.~\ref{duichen+6}.

Let $a$ and $b$ be complements of each other, i.e., $a+b=10$. Then $a=10-b$, $-a=-10+b$. So adding a number is equal to the carry minus its complement (with respect to 10), while abstracting a number is equal to the abdication minus its complement (with respect to 10). Thus, $+1$ is equal to $-9$ after carry, $-1$ is equal to $+9$ after abdication; $+2$ is equal to $-8$ after carry, and $-2 $ is equal to $+8$ after abdication, and so on.
Carry and abdication correspond to the increase and decrease of the family number of the nine-palace diagram.
With the help of the concept of complements, we can expand the basic diagram of addition and subtraction. We don't need draw all of these extended diagrams, but just summarize the following proposition.

\begin{prop}\label{jiajiantushiguanxi}
The addition diagram of each number may be aided by the subtraction diagram of its complement, in which case it carries; The subtraction diagram of each number may be aided by the addition diagram of its complement,  in which case it abdicates.
\end{prop}

The \emph{successive addition and subtraction} is to add or subtract again and again in order, that is, to add or subtract number by number. For example, in order to calculate $1 + 2 + 3 + 4 $, one can start with 1, first add $2 $, then add $3 $, and then add $4 $. What we're talking about here is sequential addition on the nine-palace diagram, that is, adding numbers sequently by the basic diagrams of  addition and subtraction for numbers.

\begin{exmp}
Calculate the following combination of addition and subtraction using a nine-palace diagram:
$1-2-9-8-7-6+8-3+5-6$.
\end{exmp}

\begin{figure}[!htbp] 
\centering
\begin{subfigure}[t]{.4\linewidth}
\centering
\begin{tikzpicture}[scale=0.6]
\draw[step=2] (0,0) grid (4,4);
\draw[fill] (-2,4) circle (3pt)  (0,4) circle (3pt)  (2,6) circle (3pt);
\coordinate [label=90:$0$](O) at (-2,4);
\coordinate [label=135:$0$](A) at  (0,4);
\coordinate [label=-45:$-1$](B) at   (2,6);
\draw[line width=3pt,->,>=latex](O)--(A);
\draw[line width=3pt,->,>=latex](A)--(B);
\end{tikzpicture}
\caption{$0+1-2$}
\label{fig:hunheyunsuan-a}
\end{subfigure}
\quad
\begin{subfigure}[t]{.4\linewidth}
\centering
\begin{tikzpicture}[scale=0.6]
\draw[step=2] (0,0) grid (4,4);
\draw[fill] (-2,4) circle (3pt)  (2,2) circle (3pt)  (2,4) circle (3pt)(4,0) circle (3pt);
\coordinate [label=45:$-1$](A) at (4,0);
\coordinate [label=90:$-1$](O) at (-2,4);
\coordinate [label=90:$-2$](B) at (2,4);
\coordinate [label=45:$-3$](C) at  (2,2);
\draw[line width=3pt,->,>=latex](A)--(O);
\draw[line width=3pt,->,>=latex](O)--(B);
\draw[line width=3pt,->,>=latex](B)--(C);
\end{tikzpicture}
\caption{$\cdots-9-8-7$}
\label{fig:hunheyunsuan-b}
\end{subfigure}
\quad
\begin{subfigure}[t]{.4\linewidth}
\centering
\begin{tikzpicture}[scale=0.6]
\draw[step=2] (0,0) grid (4,4);
\draw[fill] (0,2) circle (3pt)  (0,0) circle (3pt)(4,0) circle (3pt)(2,2) circle (3pt);
\coordinate [label=45:$-3$](A) at (2,2);
\coordinate [label=45:$-4$](B) at (4,0);
\coordinate [label=45:$-3$](C) at  (0,0);
\coordinate [label=45:$-3$](D) at (0,2);
\draw[line width=3pt,->,>=latex](A)--(B);
\draw[line width=3pt,->,>=latex](B)--(C);
\draw[line width=3pt,->,>=latex](C)--(D);
\end{tikzpicture}
\caption{$\cdots-6+8-3$}
\label{fig:hunheyunsuan-c}
\end{subfigure}
\quad
\begin{subfigure}[t]{.4\linewidth}
\centering
\begin{tikzpicture}[scale=0.6]
\draw[step=2] (0,0) grid (4,4);
\draw[fill] (0,2) circle (3pt)  (4,0) circle (3pt)  (4,4) circle (3pt);
\coordinate [label=45:$-3$](A) at (0,2);
\coordinate [label=45:$-3$](B) at  (4,0);
\coordinate [label=45:$-3$](C) at (4,4);
\draw[line width=3pt,->,>=latex](A)--(B);
\draw[line width=3pt,->,>=latex](B)--(C);
\end{tikzpicture}
\caption{$\cdots+5-6$}
\label{fig:hunheyunsuan-d}
\end{subfigure}
\caption{Calculate $1-2-9-8-7-6+8-3+5-6$ using the nine-palace diagram}
\label{fig:hunheyunsuan}
\end{figure}
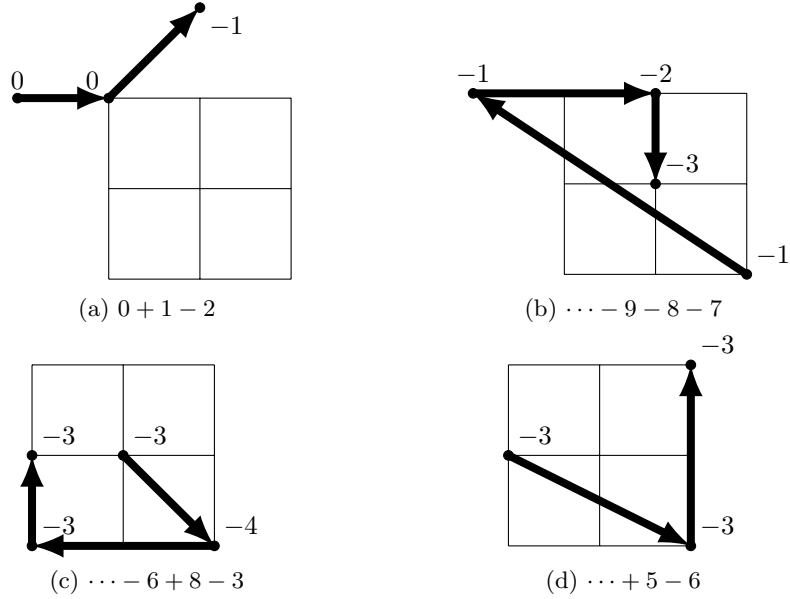

\begin{solution}
Let's start with the point $0$ in the primitive nine-palace diagram. If you want to add $1$, you move backwards one space, so you get to the upper left point. Since it is still in the $0$-family, $0$ is marked next to the dot $1$.
Subtracting $2$ is a small step of moving backwards, so we reach the point $9$ (the botton right corner) of the $1$-family diagram. Next to this point, we mark the family number  $-1$ we have currently reached.
See Fig.~\ref{fig:hunheyunsuan-a}.
\par
Subtracting $9$ reaches the position $0$, where the label is still $-1$ because the family number has not changed.
Subtracting $8$ is the same thing as subtracting $10$ and simultaneously adding $2$.  So,  the dot advances two squares to reach the top midpoint, but the family number decreases by $1$. The income family number $-2$ is used as the mark of the top midpoint.
Subtracting $7$ is the same as subtracting $10$ and simultaneously adding $3$. Thus, clicking down one grid to reach the center of the nine-palace diagram must reduce the family number by another $1$. Thus the family number $-3$ will be used as the center mark.
See Fig.~\ref{fig:hunheyunsuan-b}.
\par
Subtracting $6$ is the same as subtracting $10$ and simultaneously adding $4$. Thus, the dot goes down a little to the bottom right corner, and the family number is reduced by another $1$ to become $-4$. Use $-4$ as the mark in the lower right corner.
Adding $8$ is adding $10$ and  simultaneously subtracting $2$. Thus, the dot goes back two spaces to reach the left bottom corner, and the family number is increased by $1$ to become $-3$. Mark $-3$ at the left bottom corner.
Subtracting $3$ is equivalent to going up a little bit to the left midpoint with the same family number. The family number $-3$ is marked at the left center point.
See Fig.~\ref{fig:hunheyunsuan-b}.
\par
Adding $5$ is the same as drawing $S$ in the diagram, so when we get to the bottom right corner, the family number stays the same. Use the family number $-3$ as the mark at the bottom right corner. Subtracting $6$ is equivalent to going up two squares to reach the upper right corner with the same family number. See Fig.~\ref{fig:hunheyunsuan-d}.
\par
Finally, we reach position $3$ of family $-3$, which reads $(-3,3)=-30+3=-27$. So we have
$$1-2-9-8-7-6+8-3+5-6=-27.$$
\end{solution}

The solution process can be completed actually in the mind, and in the mental calculation process it is not needed to let the numbers $A,B,C$, etc, appear. If you can do the addition and subtraction of single-digit numbers, you can do those of multiple-digit numbers naturally.

\section{The rotation invariance property}

One of the most important rules of the nine-palace arithmetic is the so-called \emph{rotation invariance property}. What is rotation invariance? Roughly speaking, the computed bitmap is rotated by $90^\circ$ to obtain another bitmap. If a rotation of $90^\circ$ still yields a lattice, then many times of the rotation of $90^\circ$ will still yield a lattice. Because spinning three $90^\circ$ is equivalent to spinning one $90^\circ$ in the opposite direction, there is no need to distinguish the clockwise and counterclockwise directions of the rotation.
How does the zero rotate? We can expand the nine-palace diagram by four points, including adding two zero points $0$s and two ten points $10$s. When the graph rotates in the plane centered on $5$, the points $0,0,10,10$ change in a cycle successively.

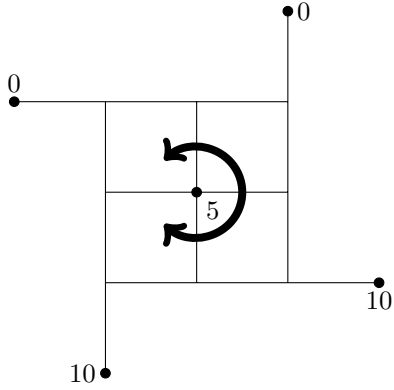
\begin{figure}[!htbp]
\centering
\begin{tikzpicture}[scale=0.6]
\draw[step=2] (0,0) grid (4,4);
\draw[fill] (-2,4) circle (3pt) (4,6) circle (3pt)(0,-2) circle (3pt) (6,0) circle (3pt)(2,2) circle (3pt);
\draw (-2,4)--(0,4)  (4,4)--(4,6) (0,0)--(0,-2)  (4,0)--(6,0);
\coordinate [label=90:$0$](a) at (-2,4);
\coordinate [label=0:$0$](a) at (4,6);
\coordinate [label=-45:$5$](a) at (2,2);
\coordinate [label=180:$10$](a) at (0,-2);
\coordinate [label=-90:$10$](a) at (6,0);
\draw[line width=3pt,<->](1.29289,1.29289) arc (-135:135:1);
\end{tikzpicture}
\caption{The expanding with 4 points and rotation of the nine-palace diagram}
\label{jiugonggexuanzhuan}
\end{figure}

Now we give a more specific and precise definition of the rotation invariance.

\begin{defn}\rm
By saying that that the nine-palace diagram has the \emph{rotation invariance property} for addition, we mean that the following property are true:  if the number $b$ is the units digit of the sum of single-digit numbers $a_1,\cdots, a_n$, then the number $b'$ is exactly the units digit of the sum of the numbers $a_1',\cdots, a_n'$, where points $a_1',\cdots, a_n',b'$ are just obtained by rotating some same times of $90^\circ$ from the original points $a_1,\cdots, a_n, b$, respectively.
\end{defn}

For example, the dot graph of $1+1=2$ represents the sum of two upper left corner points to get the top midpoint, and the dot graph of $3+3=6$ shows the sum of two upper right corner points to get the right midpoint, and the dot graph of $9+9=18$ exhibits the sum of two lower right corner points to get the lower midpoint.
The dot graph $7+7=14$ represents that two lower-left corners are added together to reach the left midpoint .

\begin{defn}\rm
By saying that the nine-palace diagram has the \emph{rotation invariance property} for subtraction, we mean that the following property is true: if the difference $a_1-\cdots -a_n$ modulo $10$ of units  digits numbers $a_1,\cdots, a_n$ gives $b$, and in the nine-palace diagram the points $a_1,\cdots, a_n, b$ are rotated simultaneously by some same times of $90^\circ$  to obtain the points $a_1 ',  \cdots, a_n ', b'$ respectively, then $b'$ is just the difference $a_1'-\cdots -a_n'$ modulo $10$.
\end{defn}

\begin{defn}\rm
Let $k$ be an integer. We say that the nine-palace diagram has the rotation invariance property if the following holds: If $b$ is the units digit of the $k$ multiple of a single digit  $a$ and if in the nine-palace diagram, the points $a,b$ are rotated simultaneously by some same times of $90^\circ$ to obtain the points $a',b'$ respectively, then the number $b'$ is exactly the units digits of the $k$ times of $a'$.
If for any integer $k$, the nine-palace diagram has the rotation invariance property for $k$ multiples, then we say that the multiples diagram has the \emph{rotation invariance property} for multiplication.
\end{defn}

\begin{figure}[!htbp] 
\centering
\begin{tikzpicture}[scale=0.8]
\draw[step=2] (0,0) grid (4,4);
\draw[fill] (0,2) circle (3pt) (4,2) circle (3pt)(2,0) circle (3pt)(2,2) circle (3pt)(2,4) circle (3pt);
\foreach \z in{(0,2),(4,2),(2,4),(2,0)} \draw[line width=3pt][->,>=latex](2,2)--\z;
\coordinate [label=90:$4-5$](a) at (1,2);
\coordinate [label=0:$2-5$](a) at (2,3);
\coordinate [label=-90:$6-5$](a) at (3,2);
\coordinate [label=-180:$8-5$](a) at (2,1);
\end{tikzpicture}
\caption{The rotation invariance property of subtraction}
\label{jianfaxuanzhuanbubianxing}
\end{figure}
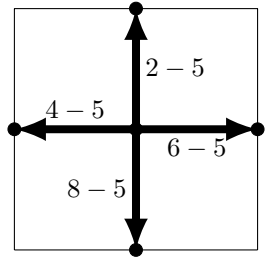

For example, as in Fig.~\ref{jianfaxuanzhuanbubianxing}, the arrow toward the right represents $6-5=1$, rotating clockwise in turn to get the downward arrow (representing $8-5=3$), the arrow toward the left (representing $4-5=-1$), the upward arrow (representing$2-5=-3$).  And if you plot these operations here, you can also see the rotation invariance of these subtraction operations.

\begin{figure}[!htbp] 
\centering
\begin{subfigure}[t]{.4\linewidth}\centering
\begin{tikzpicture}[scale=0.8]
\draw[step=2] (0,0) grid (4,4);
\draw[fill] (0,4) circle (3pt) (4,4) circle (3pt)(4,0) circle (3pt)(0,0) circle (3pt);
\draw[line width=3pt][->,>=latex](0,4)--(4,4);
\draw[line width=3pt][->,>=latex](4,4)--(4,0);
\draw[line width=3pt][->,>=latex](4,0)--(0,0);
\draw[line width=3pt][->,>=latex](0,0)--(0,4);
\coordinate [label=-180:$3\times 7$](a) at (0,2);
\coordinate [label=90:$3\times 1$](a) at (2,4);
\coordinate [label=0:$3\times 3$](a) at (4,2);
\coordinate [label=-90:$3\times 9$](a) at (2,0);
\end{tikzpicture}
\caption{The 3 times of odd numbers}
\label{beishuxuanzhuanbubianxing-a}
\end{subfigure}
\hskip40pt
\begin{subfigure}[t]{.4\linewidth}\centering
\begin{tikzpicture}[scale=0.8]
\draw[step=2] (0,0) grid (4,4);
\draw[fill] (0,2) circle (3pt) (4,2) circle (3pt)(2,0) circle (3pt)(2,4) circle (3pt);
\draw[line width=3pt][->,>=latex](0,2)--(2,4);
\draw[line width=3pt][->,>=latex](2,4)--(4,2);
\draw[line width=3pt][->,>=latex](4,2)--(2,0);
\draw[line width=3pt][->,>=latex](2,0)--(0,2);
\coordinate [label=180:$3\times 4$](a) at (0,3);
\coordinate [label=0:$3\times 6$](a) at (4,1);
\coordinate [label=-90:$3\times 8$](a) at (1,0);
\coordinate [label=90:$3\times 2$](a) at (3,4);
\end{tikzpicture}
\caption{The 3 times of of even numbers}
\label{beishuxuanzhuanbubianxing-b}
\end{subfigure}
\caption{The rotation invariance of the multiples of 3}
\label{beishuxuanzhuanbubianxing}
\end{figure}
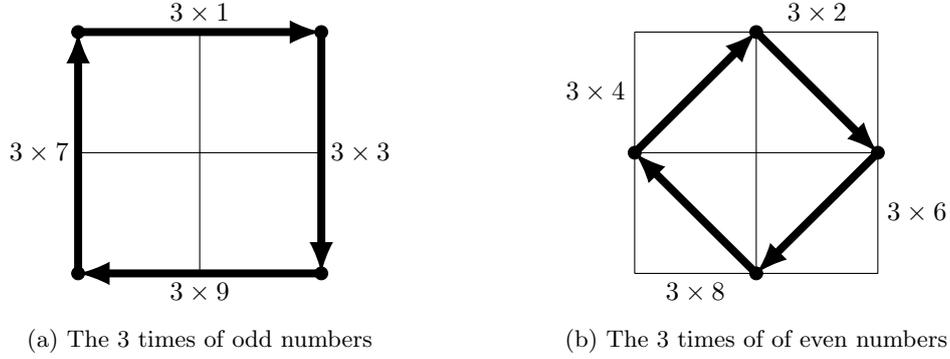

The rotation invariance of the 2 multiples  in the nine-palace diagram is easily seen from the rotation invariance of the addition of a number plus itself.
The rotation invariance of the 3 multiples can be seen in  Fig.~\ref{beishuxuanzhuanbubianxing}. In Fig.~\ref{beishuxuanzhuanbubianxing-a}, we see that the sideline with the arrow becomes the right edge line, upper edge line, and the left edge line in turn after rotation, which show that the units digits of three times of $1,3,9$ and $7$ are $3,9,7$ and $1$ respectively. Fig.~\ref{beishuxuanzhuanbubianxing-b} shows that the connecting line of two adjacent midpoints remains the connecting line of two adjacent midpoints after rotation. For example,
the units digit of three times of $2,6,8,4$ is exactly $6,8,4,2$ respectively.

\begin{thm}[The Rotation Invariance Property]
The nine-palace diagram is rotation invariant for addition, subtraction and multiplication.
\end{thm}

\begin{proof}
Let's first prove the rotation invariance of addition.
\par
Suppose that  $c$ is  the units digit of the sum $a+b$ of single-digit numbers $a$ and $b$. Let points $a, b, c$ be rotated by $90^\circ$ clockwise to obtain points $a'$, $b'$, $c'$, respectively. We want to prove that $c'$ is the the units digit of the sum $a'+b'$. For this, we discuss 10 cases where $b = 0,1,2,3,4,5,6,7,8,9$, respectively.
\par
When $b=1$: there are three arrow forms for the graph of adding 1, that is, moving rightwards one space, moving completely backwards, moving downwards one line with the right end becoming to the left end, see  Fig.~\ref{pic-plus1} and Formula~\ref{jiajiantushiguanxi}.
Let's take the example of moving one space to the right. Because the units digit of $a+1$ is $c$, in the graph, clicking $a$ to move rightwards one space produces $c$. After rotating by $90^\circ$ clockwise, moving  one space rightwards should become moving  one space downwards, and clicking $a'$ to move down one space produces $c'$. But moving down one space is the basic graph of adding 3, see~\ref{pic-plus3}. This suggests that the units digit of the sum $a'+ 3$ is $c' $. Note the point $1$ rotates by $90^\circ$ clockwise to get the point $3$, i.e., $b'=1'=3$. So, the units digit of $a'+b'$ is exactly $c'$.
\par
Similarly, we can discuss the case where $b$ is some other number. We omit the discussion of this.
\par
So this proves that the addition of two numbers is rotation invariant. Subtraction is the inverse operation of addition, so the  subtraction of any two numbers also satisfies rotation invariance property.
\par
By the induction method on the natural number $n$ we can prove that for any $n$ the corresponding addition and subtraction for $n$ have the rotation invariance property.
\par
Note that multiplication is a special case of addition since $ka=\underbrace{a+\cdots +a}_{k}$. Thus the nine-house diagram is also rotation invariant for multiplication.
\end{proof}

The rotation invariance of the nine-house diagram is important because it makes the arithmetic very simple. Roughly speaking, any numbers could be computed quickly if one can do the addition involved $1$ or $2$. We will illustrate this only with an example of addition operations.

If all involved numbers in an addition are only $5$ or $0$, then the units digit must be $0$ or $5$. Therefore, we assume that at least one of the two points used for addition is a border point. We consider, for example, $6 + 3 $, the addend here is the upper right corner point $3 $.  According to the rotation invariance, calculating $6 + 3 $ can be converted into $2 + 1 $, and that $2 + 1 = 3 $ showing the midpoint $2$ plus $1$ (moving one space rightwards)  to obtain the right corner point $3$.  The lattice $2+1=3$ rotates by $90^\circ$ clockwise  to get the lattice $6+3=9$. In other words, looking toward the right midpoint to see $6+3=9$ is the same thing as looking toward the top midpoint to see $2+1=3$. Therefore, we imagine that we are standing at the center point 5 facing the right center point, so that the corner point $3$ is like point $1$, and that $6+3$ is like point $6$ plus $1$, which shows that we move one space rightwards in the new direction. The purpose of this example is to show that adding a corner point to any point is the same as adding $1$, nothing but that the corner point would be the upper left point (the position $1$) according to the new orientation. Similarly, we can get that adding a midpoint to any point is the same thing as adding $2$, nothing but that the midpoint  would be in front of the new orientation so that it looks like the primitive $2$ in the previous orientation. Summing up, we get the following formula:

\begin{cor}[Adding by Rotation]
Given the points $a$ and $b$, we can find the units digit of their sum of $a$ and $b$ as follows.
\begin{enumerate}[(1)]
\item If $b$ is a corner point, then think of the corner point $b$ as point $1$ towards the midpoint directly adjacent to the corner point, and $a+b$ is equivalent to the original point $a$ moving one space according to the new orientation. \emph{In short, adding a corner point is like adding $1$ alaong the new direction}.
\item If $b$ is a midpoint, then look towards the midpoint and think of it as point $2$, and $a+b$ is equivalent to the original point $a$ moving two space according to the new orientation. \emph{In short, adding the midpoint is like adding $2$ along the new direction}.
\end{enumerate}
\end{cor}

The sum of two numbers is carried at most $1$ time, and whether there is a carry depends on whether the units digit of the sum is receded, which means that $1$ must be carried when the point moves forward, otherwise there will be no carry. With the above formula, we can add any number on a nine-palace diagram. Take a simple example below.

\begin{exmp}
Mentally sum up five points  $5, 3, 9, 4, 8$.
\end{exmp}

\begin{solution}
Draw five points $5, 3, 9, 4, 8$ in the primitive nine-palace diagram, see Fig.~\ref{fig:+1+2fa}.
\begin{figure}[!htbp] 
\centering
\begin{tikzpicture}[scale=0.6]
\draw[step=2] (0,0) grid (4,4);
\draw[fill] (2,2) circle (3pt) (4,4) circle (3pt) (4,0) circle (3pt) (0,2) circle (3pt)(2,0) circle (3pt);
\coordinate [label=135:$5$](B) at (2,2);
\coordinate [label=135:$3$](B) at (4,4);
\coordinate [label=45:$9$](B) at (4,0);
\coordinate [label=135:$4$](A) at  (0,2);
\coordinate [label=45:$8$](A) at  (2,0);
\end{tikzpicture}
\caption{Find the sum of five points $53948$}
\label{fig:+1+2fa}
\end{figure}
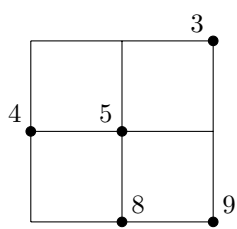
The point $5$ is the center of the nine-palace diagram. Adding $3$ is the same thing as adding $1$ when facing the right side, and so  that you get to the lower midpoint. Adding $9$, you're going to add $1$ facing downwards, and you're going to the  left bottom, and you're going to carry by 1. Adding $4$ is the same thing as adding $2$ when facing the left side, and you will get to the top left point, and you will carry by 1 again. Adding $8$ is the same thing as adding $2$ as facing downwards, and you get to the bottom corner point. Therefore, the units digit of the sum is equal to $9$. Since there are two times of carries, the sum is equal to $29$, i.e
$$5+3+9+4+8=29.$$
\end{solution}


\begin{exmp}\label{feituduichenqiuheliti3}
By the nine-palace diagram sum up $8+3+4+4+9+6+8+7+9+1.$
\end{exmp}

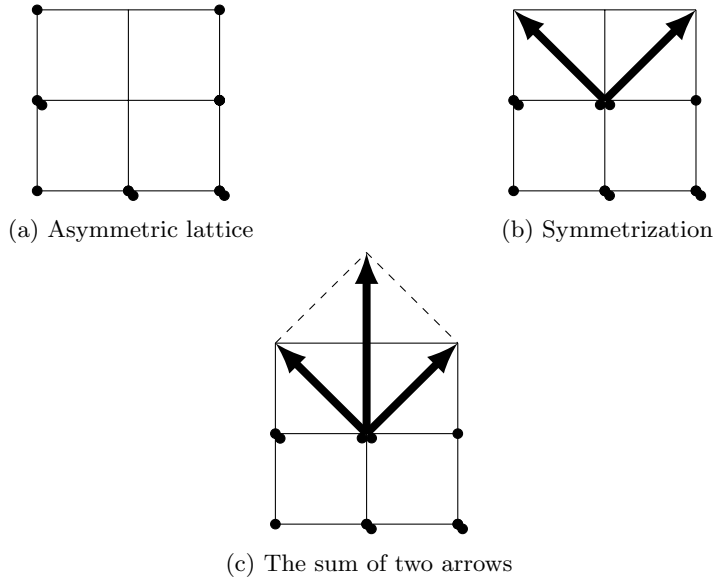
\begin{figure}[!htbp] 
\centering
\begin{subfigure}[t]{.4\linewidth}\centering
\begin{tikzpicture}[scale=0.6]
\draw[step=2] (0,0) grid (4,4);
\draw[fill] (0,0) circle (3pt)(2,0) circle (3pt) (4,0) circle (3pt)(4,2) circle (3pt)(4,2) circle (3pt)(0,4) circle (3pt)(4,4) circle (3pt);
\draw[fill] (0,2) circle (3pt) [shift={(3pt,-3pt)}] (0,2) circle (3pt);
\draw[fill] (2,0) circle (3pt) [shift={(3pt,-3pt)}] (2,0) circle (3pt);
\draw[fill] (4,0) circle (3pt) [shift={(3pt,-3pt)}] (4,0) circle (3pt);
\end{tikzpicture}
\caption{Asymmetric lattice}
\label{fig:duichenhua-qiuhe-2-a}
\end{subfigure}
\quad
\begin{subfigure}[t]{.4\linewidth}\centering
\begin{tikzpicture}[scale=0.6]
\draw[step=2] (0,0) grid (4,4);
\draw[fill] (0,0) circle (3pt)(2,0) circle (3pt) (4,0) circle (3pt)(4,2) circle (3pt);
\draw[fill] (0,2) circle (3pt) [shift={(3pt,-3pt)}] (0,2) circle (3pt);
\draw[fill] (2,0) circle (3pt) [shift={(3pt,-3pt)}] (2,0) circle (3pt);
\draw[fill] (4,0) circle (3pt) [shift={(3pt,-3pt)}] (4,0) circle (3pt);
\draw[fill] (1.9,1.9) circle (3pt) [shift={(6pt,0pt)}](1.9,1.9) circle (3pt);
\draw[line width=3pt][->,>=latex](2,2)--(0,4);
\draw[line width=3pt][->,>=latex](2,2)--(4,4);
\end{tikzpicture}
\caption{Symmetrization}
\label{fig:duichenhua-qiuhe-2-b}
\end{subfigure}
\quad
\begin{subfigure}[t]{.4\linewidth}\centering
\begin{tikzpicture}[scale=0.6]
\draw[step=2] (0,0) grid (4,4);
\draw[fill] (0,0) circle (3pt)(2,0) circle (3pt) (4,0) circle (3pt)(4,2) circle (3pt);
\draw[fill] (0,2) circle (3pt) [shift={(3pt,-3pt)}] (0,2) circle (3pt);
\draw[fill] (2,0) circle (3pt) [shift={(3pt,-3pt)}] (2,0) circle (3pt);
\draw[fill] (4,0) circle (3pt) [shift={(3pt,-3pt)}] (4,0) circle (3pt);
\draw[fill] (1.9,1.9) circle (3pt) [shift={(6pt,0pt)}](1.9,1.9) circle (3pt);
\draw[dashed] (0,4)--(2,6)--(4,4);
\draw[line width=3pt][->,>=latex](2,2)--(0,4);
\draw[line width=3pt][->,>=latex](2,2)--(4,4);
\draw[line width=3pt][->,>=latex](2,2)--(2,6);
\end{tikzpicture}
\caption{The sum of two arrows}
\label{fig:duichenhua-qiuhe-2-c}
\end{subfigure}
\caption{Find the sum by symmetry in the nine-palace diagram (example~\ref{feituduichenqiuheliti3})}
\label{fig:duichenhua-qiuhe-2}
\end{figure}

\begin{solution}
Here we're adding 10 numbers. Begin by drawing the dots $8344968791$ on the grid, as shown in the Fig.~\ref{fig:duichenhua-qiuhe-2-a}.
It is easy to see that the lattice is asymmetric. To make the lattice symmetric, simply move the points $1$ and $3$ to $5$. From point 5 draw two arrows pointing to $1$ and $3$ respectively, as shown in Fig.~\ref{fig:duichenhua-qiuhe-2-b}. The sum of these two arrows represents $-6$, see Fig.~\ref{fig:duichenhua-qiuhe-2-c}. The center of symmetry of the symmetric lattice is $67$, and the midpoint of line segment  $67$ is $6.5$,  so the sum of the symmetric lattice is $6.5$. So the sum required by the original question is $65-6=59$.
\end{solution}

\section{Finding the sum by the barycenter method}

To calculate the sum of $n$ single digits $a_1, a_2,\ldots, a_n$, the average can be calculated:
$$a=\frac{a_1+a_2+\cdots+a_n}{n},$$
and then multiplied it by $n$. If the mean is calculated in the nine-palace diagram, it is possible that the fraction $\frac{m}{n}$ cannot be represented. For this reason, we need to refine the nine grids.  The so-called \emph{$n$-refinement} is to further divide each small cell in the original nine-palace diagram into smaller cells of $n\times n$. The refinement results in many new intersections, each representing a number (possibly a fraction).
How do you know what each point represents? It can be calculated as follows:

\begin{prop}
In the $n$-refined nine-palace diagram, it increase $\frac{1}{n}$ when moving 1 small space rightwards in the horizontal direction, and $\frac{3}{n}$ when moving 1 small space downwards in the vertical direction.
\end{prop}

\begin{figure}[!htbp]
\centering
\begin{tikzpicture}[scale=0.7]
\draw[step=1] (0,0) grid (6,6);
\foreach \x in {0,3,6} \foreach \y in {0,3,6} \draw[fill] (\x,\y) circle (3pt);
\draw[line width=3pt](0,0)--(0,3)--(0,6)--(3,6)--(6,6)--(6,3)--(6,0)--(3,0)--cycle;
\draw[line width=3pt](0,3)--(3,3)--(6,3)(3,0)--(3,3)--(3,6);
\coordinate [label=-45:$1$](A) at (0,6);
\coordinate [label=-45:$2$](A) at (0,5);
\coordinate [label=-45:$3$](A) at (0,4);
\coordinate [label=-45:$4$](A) at (0,3);
\coordinate [label=-45:$5$](A) at (0,2);
\coordinate [label=-45:$6$](A) at (0,1);
\coordinate [label=-45:$7$](A) at (0,0);
\coordinate [label=-45:$2$](A) at (3,6);
\coordinate [label=-45:$3$](A) at (3,5);
\coordinate [label=-45:$4$](A) at (3,4);
\coordinate [label=-45:$5$](A) at (3,3);
\coordinate [label=-45:$6$](A) at (3,2);
\coordinate [label=-45:$7$](A) at (3,1);
\coordinate [label=-45:$8$](A) at (3,0);
\coordinate [label=-45:$3$](A) at (6,6);
\coordinate [label=-45:$4$](A) at (6,5);
\coordinate [label=-45:$5$](A) at (6,4);
\coordinate [label=-45:$6$](A) at (6,3);
\coordinate [label=-45:$7$](A) at (6,2);
\coordinate [label=-45:$8$](A) at (6,1);
\coordinate [label=-45:$9$](A) at (6,0);
\end{tikzpicture}
\caption{$3$-refined nine-palace diagram}
\label{fig:xifen-jiugongtu}
\end{figure}
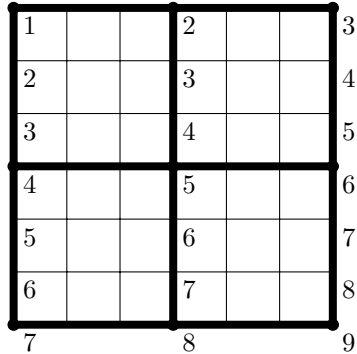

For example, Fig.~\ref{fig:xifen-jiugongtu} is a $3-$refined nine-palace diagram, in which
it increase $\frac{1}{3}$ when moving 1 small space rightwards in the horizontal direction, and $1$ when moving 1 small space downwards in the vertical direction. Points with integer values are marked, and the remaining points represent fractions.
For example, the intersection of the horizontal line $2$ and the vertical line $5$ represents  the number $3+\frac{1}{3}=\frac{4}{3}$.

Imagine that all the points on the diagram are balls of equal mass, equal weight and the same size. The center of mass of some of these points, the center of mass of all the spheres corresponding to these points, is physically the point of support that keeps the system of particles in equilibrium, and mathematically it is exactly the average of these numbers. So we can calculate the average by finding the center of mass, from which we can get the sum of those numbers. This summation method is called the \emph{barycenter method}.

First, we look at the barycenter method of summing up two points. In order to find the sum of two points $a$ and $b$ in the nine-palace diagram, we assume that the primitive nine-palace  diagram is refined by $2$ in our thinking, in which we  connect the points $a$ and $b$, and take the midpoint $c$ of the segment $ab$. Then $c$ is the center of gravity of points $a$ and $b$. Read out the corresponding number of the center of gravity $c$ in the $2$-refined nine-palace diagram, and then multiply it by $2$, so as to obtain the sum of the number $a,b$. For example, points $6$ and $7$ have the center of gravity in the point $5$ plus $\frac {1} {2} $ downwards, which represents the number $5 + \frac {3} {2} $.  As a result, we have
$$6+7=2\times 5+3=13.$$

Now we turn to the barycenter method of summing up four points.
In order to find the sum of the four points $a, b, c$ and $d$ in the nine-palace diagram, we assume in our thinking that  the primitive nine-palace  diagram is refined by $4$, in which we connect the points $a$ and $b$ and take the midpoint $p$ of the segment $ab$, connect points $c$ and $d$ and take the  midpoint $q$ of the segment $cd$, at last, connect $p$ and $q$ and take the midpoint $r$ of the segment $pq$. Then, $r$ is the center of gravity of four points $a, b, c, d$, that is, the average of the four numbers. Read out the number corresponding to the center of gravity $r$ in the $4$-refined nine-palace diagram, and then multiply it by $4$, so as to get the sum of four points $a, b, c, d$. For example, $4, 7, 8, 9 $'s center of gravity is at the point $8 $ with moving $\frac {1} {4}$ small grids leftwards and then moving $\frac {1} {4} $ small grids upwards, which represents the number $8-\frac{1}{4}-3\times \frac{1}{4}=7$. As a result,
$$4+7+8+9=4\times 7=28.$$

Sometimes it's not easy to figure out exactly where the center of gravity is in some sloping line. In order to overcome this difficulty, we use the method of obtaining the center of gravity in the vertical and horizontal directions respectively and then synthesizing the both. Specifically,
first, we move the points on each horizontal line in parallel to any vertical line, and then calculate the center of gravity of the points after moving on this vertical line. The resulting point is called the \emph{vertical center of gravity}, which represents the position of the final center of gravity in the vertical direction of the nine-palace diagram.
Second, we move in parallel the points on each vertical line to any horizontal line centrally, and then calculate the center of gravity of the points after moving on this horizontal line. The resulting point is called the \emph{horizontal center of gravity}, which represents the position of the final center of gravity in the horizontal direction of the nine-palace diagram.
Finally, the vertical center of gravity and the horizontal center of gravity are synthesized together, that is, the barycenters in two directions are comprehensively considered.

The vertical center of gravity and the horizontal center of gravity are calculated on the same line respectively. Let's take the horizontal center of gravity as an example. For example, there are three $7$s and five $8$s, which are located in two adjacent positions on the third horizontal line of the nine-palace diagram. Since the ratio of the weights of the points contained in the two positions is $3:5$, the ratio of the distances between the center of gravity and the two positions is equal to $5:3$, that is, the center of gravity is at the point $8$ moving $\frac{3}{8}$ of a grid leftwards.

When there are points on a line that are not adjacent to each other, we can use the advance and retreat method to bring them adjacent. For example, the dot series  $77899999$ is on the third horizontal line of the nine-palace diagram, but $7$ is not adjacent to $9$. Note that $77899999= 88888999$. Thus we can eliminate two $7$s and two $9$s, and add four $8$s to get five $8$s and three $9$s. Now it is contiguous, thus it is easy to figure out that the center of gravity, which is the point $8$ with moving $\frac{3}{8}$ of a gird rightwards.

To sum up, we get the method of finding the barycenter of some points on a line: (1) if the points are non-adjacent, then move the outer points to the middle position in pairs to eliminate the non-adjacent points; (2) if the points are all adjacent, then the center of gravity is close to the point with more overlap, and the ratio of proximity is the numbers of the points with less overlap divided by the total point numbers.
The following example shows that the barycenter  method is a very effective way to find the sum of a large amount of numbers.
\begin{exmp}\label{exp:zhongxinfa}
By  the barycenter  method sum up $1, 4, 4, 7, 7, 7, 2, 5, 8, 8$.
\end{exmp}

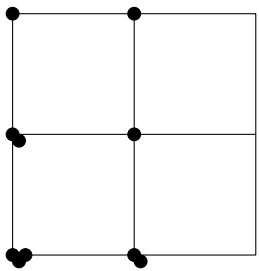
\begin{figure}[!htbp] 
\centering
\begin{tikzpicture}[scale=0.8]
\draw[step=2] (0,0) grid (4,4);
\draw[fill] (0,4) circle (3pt) (2,2) circle (3pt)(2,4) circle (3pt);
\draw[fill] (0,2) circle (3pt) [shift={(3pt,-3pt)}] (0,2) circle (3pt);
\draw[fill] (2,0) circle (3pt) [shift={(3pt,-3pt)}] (2,0) circle (3pt);
\draw[fill] (0,0) circle (3pt) [shift={(3pt,-3pt)}] (0,0) circle (3pt)[shift={(3pt,3pt)}] (0,0) circle (3pt);
\end{tikzpicture}
\caption{Addition by the barycenter method (Example~\ref{exp:zhongxinfa})}
\label{liti-zhongxinfa-qiuhe}
\end{figure}

\begin{solution}
First, draw the ten numbers on the nine-palace diagram, as shown in  Fig.~\ref{liti-zhongxinfa-qiuhe}. Concentrate the 10 points on the horizontal line, then we get just  $6$ left  points,  $6$ middle  points, so the horizontal center of gravity is at left moving $\frac{4}{10}$ steps rightwards. The ten points were concentrated on the vertical line, getting that the numbers of the upper, middle and lower points are $2, 3$ and $5 $ respectively. When the points are concentrated towards the middle point, the $7$  middle points and the $3$ lower  points are obtained. Therefore, the vertical center of gravity is lower $\frac{3}{10}$ steps then the middle point.
So the total center of gravity is at left middle $4$ moving $\frac{4}{10}$ steps rightwards and simultaneously moving $\frac{3}{10}$ steps downwards, which means the average numbers is $4+\frac{4}{10}+3\times \frac{3}{10}=5.3$. Therefore, the sum we want to find is $10\times 5.3=53.$
\end{solution}

Let us look at an other example.
\begin{exmp}\label{exmp:chengjihe}
Calculate the sum of the products of numbers with the nine-palace diagram:
$$2\times 6+4\times 9+3\times 8+4\times 3+7\times 5.$$
\end{exmp}

\begin{solution}
First, point out the following twenty points on the nine-palace diagram: $2$ times of the point $6$, $4$ times of the point $9$, $3$ times of the point $8$, $4$ times of the point  $3$, $7$ times of the point $5$.
The appearing frequency of the nine points from $1$ to $9$ is $0, 0, 4, 0, 7, 2, 0, 3, 4$ respectively.

\begin{table}[!htbp]
\centering
\begin{subfigure}[t]{.8\linewidth}
\centering
\begin{tabular}{ccc|ccc}
   0 & 0 & 4& 4& 0 &\textrm{up}\\
   0 & 7 & 2& 9& 17 &\textrm{middle}\\
   0 & 3 & 4& 7& 3 &\textrm{dowm}\\
   \hline
   0 & 10 & 10& \\
    \textrm{left} & \textrm{middle} & \textrm{right}& \\
       \hline
\end{tabular}
\end{subfigure}
\quad
\begin{subfigure}[t]{.9\linewidth}
\centering
\parbox[b]{0.8\linewidth}{\emph{Horizontal barycenter:  the middle with $\frac{10}{20}$ rightwards;  \\Vertical barycenter:  the middle with $\frac{3}{20}$  downwards;}}
\parbox[b]{0.8\linewidth}{\emph{Total barycenter: the center with $\frac{10}{20}$ rightwards and $\frac{3}{20}$  downwards.}}
\end{subfigure}
\caption{Addition by the barycenter method (Example~\ref{exmp:chengjihe})}
\label{table:chengjihe}
\end{table}

As Table~\ref{table:chengjihe} shows,
after the 20 points are concentrated on the horizontal line, the point numbers' distribution is $0, 10, 10$, so the horizontal center of gravity is the middle point moving $\frac{10}{20}$ step rightwards. When these 20 points are concentrated on the vertical line, the point numbers' distribution is $4, 9, 7$, and when they are concentrated on the midpoint, the point numbers' distribution is $0, 17, 3$, so the vertical center of gravity is the middle point moving $\frac{3}{20}$ step downwards.
So  the center of gravity of the 20 points is  the center point $5 $ moving $\frac{10}{20}$ step rightwards and simultaneously moving $\frac{3}{20}$ steps downwards, which means that the average of 20 numbers is $5+\frac{10}{20}+3\times \frac{3}{20}=5\frac{19}{20}$. So the required sum of the products is $20\times 5\frac{19}{20}=20\times 5+19=119.$
\end{solution}

\section{The dot matrix addition}

Representing two units digit numbers on a nine-palace diagram and finding the units digit of their sum and product is equivalent to addition and multiplication of modulus 10 \cite{chenjingrun-1,chenjingrun-2}. The addition and multiplication of the remaining classes of all integers modulo $10$ form semigroups respectively. For some discussion of semigroups, see for example \cite{Zhu-7,Zhu-9,Zhu-10,Zhu-1,Zhu-3,Zhu-4,Zhu-6}.
The array of points representing some units digits and the units digits of their sum on a nine-palace diagram is called the (addition) \emph{dot matrix}. The study of the dot matrices with regularity is the basis of quick addition calculation by using the nine-palace diagram.
There are indeed many beautiful dot matrices, and we will present only a few of them.

\begin{prop}
When summing up two center points or two opposite points, it has carry $1$ directly, while the units digit goes back to zero.
\end{prop}

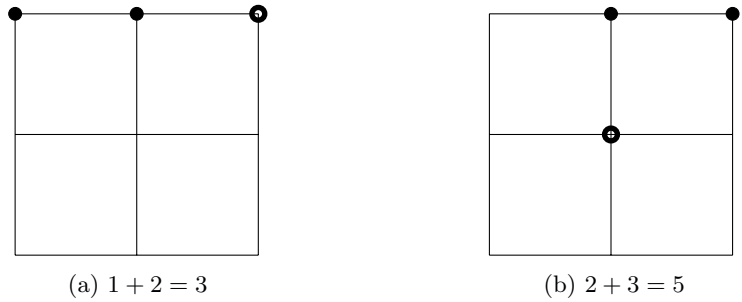
\begin{figure}[!htbp] 
\centering
\begin{subfigure}[t]{.4\linewidth}\centering
\begin{tikzpicture}[scale=0.8]
\draw[step=2] (0,0) grid (4,4);
\draw[fill] (0,4) circle (3pt) (2,4) circle (3pt);
\draw[line width=2pt] (4,4) circle (3pt);
\end{tikzpicture}
\caption{$1+2=3$}\label{fig:jiaodian+zhongdian-dianzhentu-a}
\end{subfigure}
\quad
\begin{subfigure}[t]{.4\linewidth}\centering
\begin{tikzpicture}[scale=0.8]
\draw[step=2] (0,0) grid (4,4);
\draw[fill] (2,4) circle (3pt)(4,4) circle (3pt);
\draw[line width=2pt] (2,2) circle (3pt);
\end{tikzpicture}
\caption{$2+3=5$}\label{fig:jiaodian+zhongdian-dianzhentu-b}
\end{subfigure}
\caption{The sum of corner points and adjacent midpoints (Two dot matrices)}
\label{fig:jiaodian+zhongdian-dianzhentu}
\end{figure}

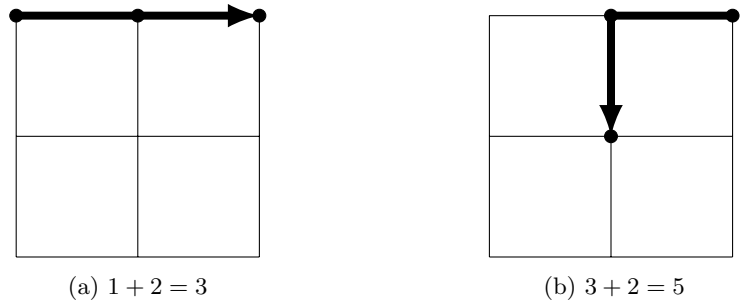
\begin{figure}[!htbp] 
\begin{center}
\begin{subfigure}[t]{.4\linewidth}\centering
\begin{tikzpicture}[scale=0.8]
\draw[step=2] (0,0) grid (4,4);
\draw[fill] (0,4) circle (3pt) (2,4) circle (3pt)(4,4) circle (3pt);
\draw[line width=3pt](0,4)--(2,4)[->,>=latex]--(4,4);
\end{tikzpicture}
\caption{$1+2=3$}\label{fig:jiaodian+zhongdian-a}
\end{subfigure}
\quad
\begin{subfigure}[t]{.4\linewidth}\centering
\begin{tikzpicture}[scale=0.8]
\draw[step=2] (0,0) grid (4,4);
\draw[fill] (2,2) circle (3pt) (2,4) circle (3pt)(4,4) circle (3pt);
\draw[line width=3pt](4,4)--(2,4)[->,>=latex]--(2,2);
\end{tikzpicture}
\caption{$3+2=5$}\label{fig:jiaodian+zhongdian-b}
\end{subfigure}\centering
\caption{The arrow points to the sum (Improved dot matrices)}
\label{fig:jiaodian+zhongdian}
\end{center}
\end{figure}

When adding two numbers, two black dots can be replaced by a white dot. But for the sake of visualization and memorization, the dots can be connected  appropriately to get a better picture with the white dots becoming black dots and with an arrow pointing to it. If connecting the dots of the graph of $1+2=3$, one gets a line, while the dot graph of $3+2=5$ looks like a small carpenter's square with an arrow pointing to the sum $5$, as shown in Fig.~\ref{fig:jiaodian+zhongdian}.

\begin{prop}\label{prop:2chongdianhe} 
In the nine-palace diagram, running clockwise in the border square, (1) any corner point plus the adjacent midpoint gets the next corner (or briefly: corner-middle-corner on a line); (2) the middle point plus the adjacent corner points equals to the center of the nine-palace diagram (or: middle-corner-center along a broking line); (3) there is no carry when the point goes forward and it will carry $1$ when the point goes backward.
\end{prop}

\begin{figure}[!htbp] 
\centering
\begin{subfigure}[t]{.4\linewidth}\centering
\begin{tikzpicture}[scale=0.8]
\draw[step=2] (0,0) grid (4,4);
\draw[fill] (2,4) circle (3pt) (0,2) circle (3pt)(4,2) circle (3pt);
\draw[line width=3pt](0,2)--(2,4)[->,>=latex]--(4,2);
\end{tikzpicture}
\caption{$4+2=6$}\label{fig:zhongzhongzhong-a}
\end{subfigure}
\quad
\begin{subfigure}[t]{.4\linewidth}\centering
\begin{tikzpicture}[scale=0.8]
\draw[step=2] (0,0) grid (4,4);
\draw[fill] (2,4) circle (3pt) (4,2) circle (3pt)(2,0) circle (3pt);
\draw[line width=3pt](2,4)--(4,2)[->,>=latex]--(2,0);
\end{tikzpicture}
\caption{$2+6=8$}\label{fig:zhongzhongzhong-b}
\end{subfigure}
\caption{The sum of two consecutive midpoints yields the next midpoint}
\label{fig:zhongzhongzhong}
\end{figure}
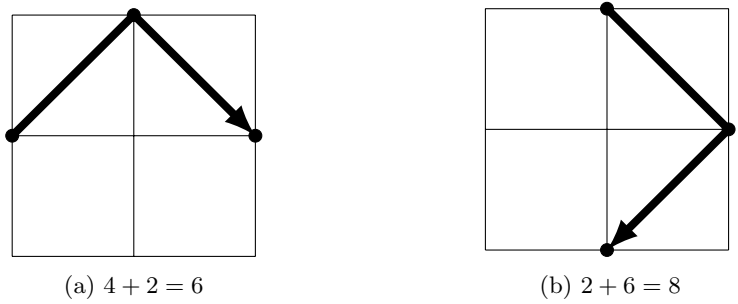

Now consider how the nine-palace diagram shows when summing up two consecutive midpoints such as $4$ and $2$. Because $4+2=6$, we see that $4, 2$ and $6$ happen to be the three midpoints on the border square in a clockwise fashion. According to the addition rotation invariance of the nine-palace diagram, the units digit of the sum of any two consecutive midpoints is equal to the next adjacent midpoint. The rounding rule is the same as before: no carry when going forward,  it will carry $1$ when going backwards.
Fig.~\ref{fig:zhongzhongzhong} shows the rule of summing up two consecutive corner points. We put the conclusion into a formula:

\begin{prop} 
The sum of two consecutive midpoints gives the next adjacent midpoint. There is no carry when going forward, and it will carry $1$ when going backwards. In short: middle-middle-middle in the square.
\end{prop}

In order to obtain some fixed diagrams able used for calculation, we often use the rotation invariance of the nine-palace diagram.

\begin{figure}[!htbp] 
\centering
\begin{tikzpicture}[scale=0.8]
\draw[step=2] (0,0) grid (4,4);
\draw[fill] (0,4) circle (3pt) (2,4) circle (3pt)(4,4) circle (3pt)(4,2) circle (3pt);
\draw[line width=3pt][->,>=latex](0,4)--(2,4)--(4,4)--(4,2);
\end{tikzpicture}
\caption{Sum up three points on one side}
\label{fig:yibianqiuhe}
\end{figure}
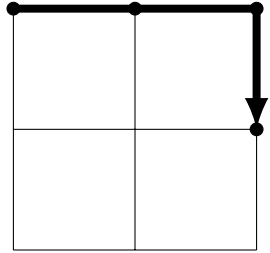

Since $1+2+3=6$, the sum of three different points on the same edge of the outer frame of the nine-palace diagram is equal to the adjacent midpoint according to the addition rotation invariance of the nine-palace diagram. Connect four dots $1236$ in sequence to get a figure similar to the head of a Chinese word which means buying, as shown in  Fig.~\ref{fig:yibianqiuhe}.  So we obtain the following formula:

 \begin{prop}
The units digit of the sum of three different points on the same edge of the outer frame is equal to the adjacent midpoint  in the nine-palace diagram. Carry according to the midpoint among the three points. In short, corner-middle-corner-middle.
 \end{prop}

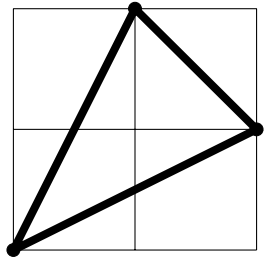
\begin{figure}[!htbp] 
\centering
\begin{tikzpicture}[scale=0.8]
\draw[step=2] (0,0) grid (4,4);
\draw[fill] (2,4) circle (3pt) (4,2) circle (3pt)(0,0) circle (3pt);
\draw[line width=3pt](2,4)--(4,2)--(0,0)(0,0)--(2,4);
\end{tikzpicture}
\caption{Three points of the permutation type such as three mid-mid-far points}
\label{zhihuanxing3dianqiuhe}
\end{figure}

Three points that are neither in the same row nor in the same column (excluding the zero) are called \emph{three points of the permutation type}. It includes two types: the diagonal three points,the three mid-mid-far points. The  so-called \emph{three diagonal points} refer to three different points on the same diagonal line of the nine-palace diagram, namely three points $159$ and three points $357$.
Any three diagonal points contain a pair of opposing points and the center of the nine-palace diagram. The sum of a pair of opposing points is $10$, and the center represents $5$. Therefore, the sum of three diagonal points is equal to $15$. The \emph{three mid-mid-far points} are actually the two consecutive midpoints on the border of the nine-palace diagram and a corner point that is not adjacent to them (in fact, it is also the corner farthest from the two midpoints). For example, $267$ constitute three mid-mid-far points, see  Fig.~\ref{zhihuanxing3dianqiuhe}. Using the advance and retreat method, it is easy to convert two adjacent midpoints into points on the diagonal, so the three mid-mid-far points can be converted into three points on the diagonal, such as $2+6+7=1+5+7$.
Therefore, the sum of the three mid-mid-far points is also equal to $15$. So we have

\begin{prop}\label{prop:zhihuanxingsandianhe}
The sum of three points of the permutation type is $15$.
\end{prop}

As shown in  Fig.~\ref{zhihuanxing3dianqiuhe},  we put the concept of three points of the permutation type and the formula of sum together with the nine-palace diagram, so that it is easy to compare and understand.

The nine non-zero points in the nine-palace diagram can be seen as being made up of three sets of  three points of the permutation type (e.g., $159$, $267$ and $348$),  so their sum is $15\times 3= $45. If a set of three sets is removed, the remaining six points, called  the \emph{six points of the permutation type}, the sum of which is $15\times 2=30$. So we get

\begin{cor}
The sum of six points of the permutation type is $30$.
\end{cor}

\begin{figure}[!htbp] 
\centering
\begin{tikzpicture}[scale=0.8]
\draw[step=2] (0,0) grid (4,4);
\draw[fill] (0,0) circle (3pt)(0,4) circle (3pt) (2,0) circle (3pt)(2,2) circle (3pt)(4,2) circle (3pt)(4,4) circle (3pt);
\end{tikzpicture}
\caption{Six points of the permutation type}
\label{fig:9-3point}
\end{figure}
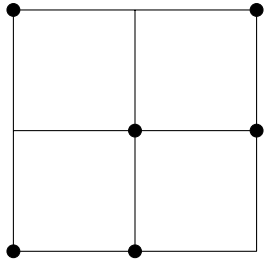

For example,  Fig.~\ref{fig:9-3point} shows that the sum of six points $135678$ is $30$.

\begin{figure}[!htbp] 
\centering
\begin{tikzpicture}[scale=0.8]
\draw[step=2] (0,0) grid (4,4);
\draw[fill] (0,2) circle (3pt)(0,4) circle (3pt) (2,4) circle (3pt)(4,4) circle (3pt);
\draw[line width=3pt][->,>=latex](0,2)--(0,4)--(2,4)--(4,4);
\end{tikzpicture}
\caption{Sum of four points of the $1234$ type}
\label{fig:1234he}
\end{figure}
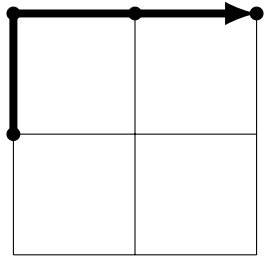

Note that $1+2+3+4=10$. When the four dots $4123$ are connected in sequence, an Arabic numeral '$7$' is formed on the border of the nine-palace diagram in a clockwise direction, as shown in  Fig.~\ref{fig:1234he}.
According to the rotation invariance, the units digit of the sum of $4123$, $2369$, $6987$, and $8741$ that form the glyph "$7$" is $0$ respectively. These four dot matrix graph are all called that of the \emph{the $1234$ type}.
According to the position of the turning point of "$7$", rotating one circle clockwise, the sum of four points of the $1234$ type is equal to $10$, $20$, $30$, $20$ and $10$ in turn. That is, when the turning point is at the corner point $1$, the carry is $1$; at the corner point $3$, carry is $2$; at the corner point $9$, the carry is $3$; at the corner point $7$, the carry is $2$. That is to say, the four point of the $1234$ type, according to the turning point, the carry rule is $1223$.
We have written these results in the following formula.

\begin{prop}\label{prop:1234he}
The sum of the four consecutive points of the mid-corner-mid-corner in the direction of the "$7$" shape is some times of $10$, and according to the turning point (or the second midpoint) the carrying law is $1223$.
\end{prop}

Eight different points on the border frame of the nine-palace diagram are connected to form two separativee "$7$" figures (e.g., $4123$ and $6987$).
So, according to Prop.~\ref{prop:1234he}, the sum has the units digit $0$, and the carry is $1+3=2+2=4$.
That is, the sum of the eight points is $10+30=20+20=40$. Or, to put it another way, the eight points are actually made up of four pairs of opposites, so the sum is $40$. We also write this consequence as a formula:

\begin{prop}
The sum of the eight points on the border frame of the nine-palace diagram is $40$.
\end{prop}

\section{A counting method of finding multiples of points}

From now on, we study the multiplication on the nine-palace diagram.
The units digits of the product of two units digits are well represented on the grid. Because of the rotation invariance for multiples, as long as we know the law of various multiples of $1$ and $2$, we can know the law of various multiples of all the numbers except $0$ and $5$. Combining these rules, we see that we can get multiples of any border point just by naturally counting numbers 1, 2, and so on. So we call this method the \emph{counting method}.

Let's first assume that the number $a$ is an odd number. Then $a=1,3,5,7,9$, that is, point $a$ is the center or corner point of the grid.

The case of the center point is particularly simple. Since $5$ multiplied by any even number is equal to a number with the units digit $0$, and $5$ multiplied by any odd number is a number whose units digits is still $5$, we get the following simple formula:

\begin{prop}\label{prop:5timesgewei}
The center point times the corner point is the center point, and the center point times a  midpoint returns to zero.
\end{prop}

Let's assume that $a=1,3,7,9$, e.g., $a$ is a corner point.

Giving a non-zero number $b $. Note that $1 \times b = b$.  If we draw a arrow from the point $1$ to $b$, then the arrow points to the units digit of $1\times b$. According to the rotation invariance of the nine-palace diagram for multiplication, one may rotate the nine-palace diagram. When the starting point of the arrow rotates to the point $a$, the end point of the arrow presents the units digit of $a\times b$. So we get

\begin{thm}[The method of finding the units digit of the multiples of a corner point]
For any corner point $a$ and any non-zero point $b$, one makes an arrow pointing to the point $b$ from the point $1$, and then rotate the arrow. When the starting point of the arrow rotates to the point $a$, the end of the arrow represents the units digit of $a\times b$.
\end{thm}
The point corresponding to the units digit of a product can be called the point of the product.

For any corner point $a$ and non-zero point $b$, the position of $a\times b$ relative to the point $a$ is the position of $b$ relative to $1$. So if we read the dot $a$ as $1$ and recount numbers in the new orientation but in the same order as the primitive nine-palace diagrams, the dot $a\times b$ is pronounced exactly as $b$. In this way, calculating $a\times b$ is like counting numbers on the diagram: If you count from $1$ to $9$, you get all the different multiples of $a$.

\begin{figure}[!htbp] 
\centering
\begin{subfigure}[t]{.4\linewidth}
\centering
\begin{tikzpicture}[scale=0.7]
\node(1) at (0,4) [circle,draw]{1};
\node(2) at (2,4) [circle,draw]{2};
\node(3) at (4,4) [circle,draw]{3};
\node(4) at (0,2) [circle,draw]{4};
\node(5) at (2,2) [circle,draw]{5};
\node(6) at (4,2) [circle,draw]{6};
\node(7) at (0,0) [circle,draw]{7};
\node(8) at (2,0) [circle,draw]{8};
\node(9) at (4,0) [circle,draw]{9};
\draw (1)--(2)--(3);
\draw (4)--(5)--(6);
\draw (7)--(8)--(9);
\draw (1)--(4)--(7) (2)--(5)--(8) (3)--(6)--(9);
\end{tikzpicture}
\caption{The multiples of $1$}
\label{fig:1beishushu}
\end{subfigure}
\quad
\begin{subfigure}[t]{.4\linewidth}
\centering
\begin{tikzpicture}[scale=0.7]
\node(1) at (0,4) [circle,draw]{7};
\node(2) at (2,4) [circle,draw]{4};
\node(3) at (4,4) [circle,draw]{1};
\node(4) at (0,2) [circle,draw]{8};
\node(5) at (2,2) [circle,draw]{5};
\node(6) at (4,2) [circle,draw]{2};
\node(7) at (0,0) [circle,draw]{9};
\node(8) at (2,0) [circle,draw]{6};
\node(9) at (4,0) [circle,draw]{3};
\draw (1)--(2)--(3);
\draw (4)--(5)--(6);
\draw (7)--(8)--(9);
\draw (1)--(4)--(7) (2)--(5)--(8) (3)--(6)--(9);
\end{tikzpicture}
\caption{The multiples of $3$}
\label{fig:3beishushu}
\end{subfigure}
\caption{The multiples of $1$ and $3$ on the nine-palace diagram}
\label{fig:13beishushu}
\end{figure}
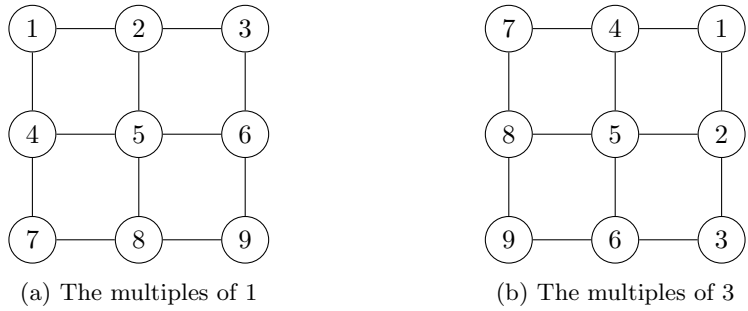

For example,  the  $3$ times of from 1 to $9$  are starting with the point $3$ followed by $6, 9,2,5,8,1,4,7$ in turn, being counting as $1$, $2$, $3$, $4$, $5$, $6$, $7$, $8$ and $9$. It is also like turning the original numbers of the primitive nine-palace diagram around $90^\circ$, as shown in  Fig.~\ref{fig:13beishushu}.

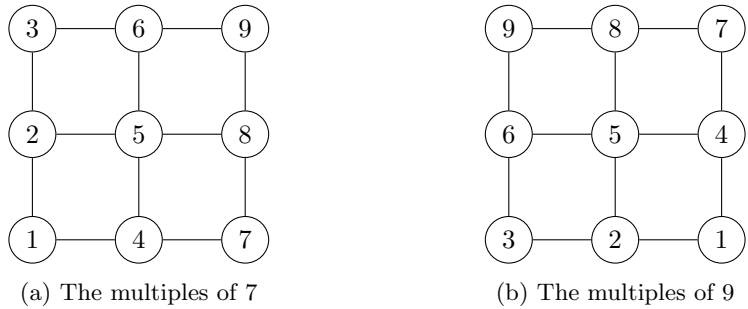
\begin{figure}[!htbp] 
\centering
\begin{subfigure}[t]{.4\linewidth}
\centering
\begin{tikzpicture}[scale=0.7]
\node(1) at (0,4) [circle,draw]{3};
\node(2) at (2,4) [circle,draw]{6};
\node(3) at (4,4) [circle,draw]{9};
\node(4) at (0,2) [circle,draw]{2};
\node(5) at (2,2) [circle,draw]{5};
\node(6) at (4,2) [circle,draw]{8};
\node(7) at (0,0) [circle,draw]{1};
\node(8) at (2,0) [circle,draw]{4};
\node(9) at (4,0) [circle,draw]{7};
\draw (1)--(2)--(3);
\draw (4)--(5)--(6);
\draw (7)--(8)--(9);
\draw (1)--(4)--(7) (2)--(5)--(8) (3)--(6)--(9);
\end{tikzpicture}
\caption{The multiples of $7$}
\label{fig:7beishushu}
\end{subfigure}
\quad
\begin{subfigure}[t]{.4\linewidth}
\centering
\begin{tikzpicture}[scale=0.7]
\node(1) at (0,4) [circle,draw]{9};
\node(2) at (2,4) [circle,draw]{8};
\node(3) at (4,4) [circle,draw]{7};
\node(4) at (0,2) [circle,draw]{6};
\node(5) at (2,2) [circle,draw]{5};
\node(6) at (4,2) [circle,draw]{4};
\node(7) at (0,0) [circle,draw]{3};
\node(8) at (2,0) [circle,draw]{2};
\node(9) at (4,0) [circle,draw]{1};
\draw (1)--(2)--(3);
\draw (4)--(5)--(6);
\draw (7)--(8)--(9);
\draw (1)--(4)--(7) (2)--(5)--(8) (3)--(6)--(9);
\end{tikzpicture}
\caption{The multiples of $9$}\label{fig:9beishushu}
\end{subfigure}
\caption{The multiples of $7$ and $9$ on the nine-palace diagram}
\label{fig:79beishushu}
\end{figure}

Multiples of $1$ are like counting numbers facing up, and multiples of $3$ are like counting numbers facing the right. Similarly, finding multiples of $7$ and $9$ would be like counting numbers facing the left and down, respectively, as shown in  Fig.~\ref{fig:79beishushu}. Note that  the number marked in the figure represents multiples, the dot with the label $1$ represent the multiplicand, and the numbers represented by each dot in the primitive nine-palace diagram represents the units digit of the product. As you see, the units digit of the product can be read directly by counting numbers from the graph.

\begin{thm}[The counting method to find the multiples of a corner point]\label{koujue:shushu1}
Finding the units digit of the multiples of a corner point is equivalent to starting with this point and counting  $1 $, $2 $, $3 $, $4 $, $5 $, $6$, $7 $, $8$ and $9 $  in turn. When we count $k$, the point we have arrived represents  the units digit of the $k$ times of the given corner point.
\end{thm}

Now let's look at the units number of the various multiples of a midpoint on the nine-palace diagram.

First, we look at multiples of $2$. Since $2\times 1=1\times 2=2$, and $2\times 2=4$, multiples of $2$ will make all corner points becomes midpoints in clockwise direction and all midpoints becomes midpoints in counterclockwise direction according to the rotation invariance of the nine-palace diagram for multiples.
Therefore, the multiples of $2$ make the points $1,2,3,4,5,6,7,8,9$ in order to become the top midpoint, the left midpoint, the right midpoint, the lower midpoint, zero, the top midpoint, the left midpoint, the right midpoint and the lower midpoint.

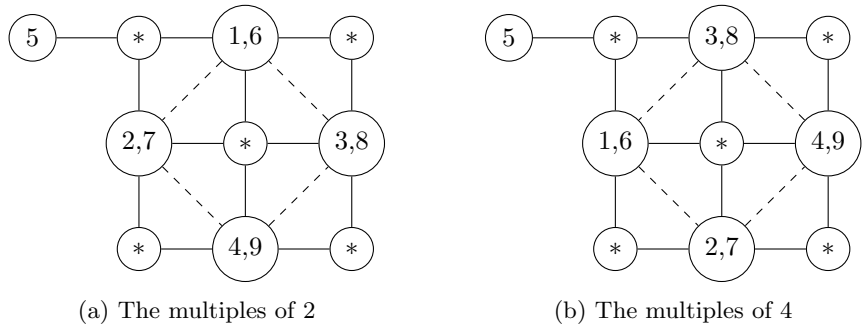
\begin{figure}[!htbp] 
\centering
\begin{subfigure}[t]{.4\linewidth}
\centering
\begin{tikzpicture}[scale=0.7]
\node(0) at (-2,4) [circle,draw]{$5$};
\node(1) at (0,4) [circle,draw]{$\ast$};
\node(2) at (2,4) [circle,draw]{1,6};
\node(3) at (4,4) [circle,draw]{$\ast$};
\node(4) at (0,2) [circle,draw]{2,7};
\node(5) at (2,2) [circle,draw]{$\ast$};
\node(6) at (4,2) [circle,draw]{3,8};
\node(7) at (0,0) [circle,draw]{$\ast$};
\node(8) at (2,0) [circle,draw]{4,9};
\node(9) at (4,0) [circle,draw]{$\ast$};
\draw (0)--(1)--(2)--(3);
\draw (4)--(5)--(6);
\draw (7)--(8)--(9);
\draw (1)--(4)--(7) (2)--(5)--(8) (3)--(6)--(9);
\draw[dashed] (2)--(4)--(8)--(6)--(2);
\end{tikzpicture}
\caption{The multiples of $2$}
\label{fig:2beishushu}
\end{subfigure}
\quad
\begin{subfigure}[t]{.4\linewidth}
\centering
\begin{tikzpicture}[scale=0.7]
\node(0) at (-2,4) [circle,draw]{$5$};
\node(1) at (0,4) [circle,draw]{$\ast$};
\node(2) at (2,4) [circle,draw]{3,8};
\node(3) at (4,4) [circle,draw]{$\ast$};
\node(4) at (0,2) [circle,draw]{1,6};
\node(5) at (2,2) [circle,draw]{$\ast$};
\node(6) at (4,2) [circle,draw]{4,9};
\node(7) at (0,0) [circle,draw]{$\ast$};
\node(8) at (2,0) [circle,draw]{2,7};
\node(9) at (4,0) [circle,draw]{$\ast$};
\draw (0)--(1)--(2)--(3);
\draw (4)--(5)--(6);
\draw (7)--(8)--(9);
\draw (1)--(4)--(7) (2)--(5)--(8) (3)--(6)--(9);
\draw[dashed] (2)--(4)--(8)--(6)--(2);
\end{tikzpicture}
\caption{The multiples of $4$}
\label{fig:4beishushu}
\end{subfigure}
\caption{The multiples of $2$ and $4$ on the nine-palace diagram}
\label{fig:24beishushu}
\end{figure}

Because even multiples of the center point return to zero, we omit the multiples of $5$ temporarily.
The $1,2,3,4$ multiples of $2$ are just on the slopping square shape, and when viewed from the center of the diagram facing the point $2$, that is, from the top of the diagram, they happen to be the front, the left, the right and the back. The $6,7,8,9$ times of $2$ are also similar to $1,2,3,4$ times of $2$, see Pic.~\ref{fig:2beishushu}. Simply saying, the multiples of $2$ are counted facing $2$, in the order of 'front, left, right and back'.  Please note: it is assumed that the person is standing at the center of the nine-palace diagram and facing the point $2$ in the original nine-palace diagram.

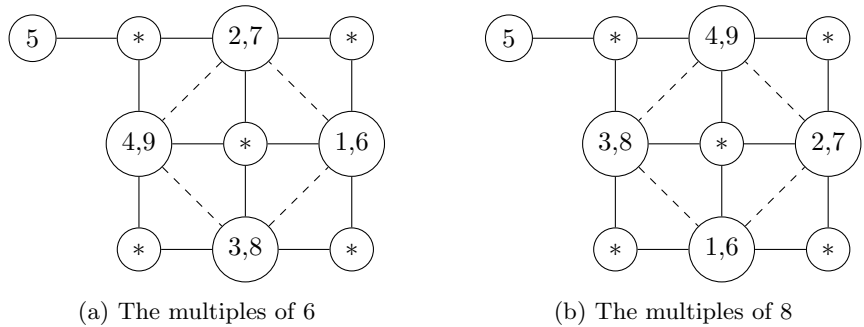
\begin{figure}[!htbp] 
\centering
\begin{subfigure}[t]{.4\linewidth}
\centering
\begin{tikzpicture}[scale=0.7]
\node(0) at (-2,4) [circle,draw]{$5$};
\node(1) at (0,4) [circle,draw]{$\ast$};
\node(2) at (2,4) [circle,draw]{2,7};
\node(3) at (4,4) [circle,draw]{$\ast$};
\node(4) at (0,2) [circle,draw]{4,9};
\node(5) at (2,2) [circle,draw]{$\ast$};
\node(6) at (4,2) [circle,draw]{1,6};
\node(7) at (0,0) [circle,draw]{$\ast$};
\node(8) at (2,0) [circle,draw]{3,8};
\node(9) at (4,0) [circle,draw]{$\ast$};
\draw (0)--(1)--(2)--(3);
\draw (4)--(5)--(6);
\draw (7)--(8)--(9);
\draw (1)--(4)--(7) (2)--(5)--(8) (3)--(6)--(9);
\draw[dashed] (2)--(4)--(8)--(6)--(2);
\end{tikzpicture}
\caption{The multiples of $6$}
\label{fig:6beishushu}
\end{subfigure}
\quad
\begin{subfigure}[t]{.4\linewidth}
\centering
\begin{tikzpicture}[scale=0.7]
\node(0) at (-2,4) [circle,draw]{$5$};
\node(1) at (0,4) [circle,draw]{$\ast$};
\node(2) at (2,4) [circle,draw]{4,9};
\node(3) at (4,4) [circle,draw]{$\ast$};
\node(4) at (0,2) [circle,draw]{3,8};
\node(5) at (2,2) [circle,draw]{$\ast$};
\node(6) at (4,2) [circle,draw]{2,7};
\node(7) at (0,0) [circle,draw]{$\ast$};
\node(8) at (2,0) [circle,draw]{1,6};
\node(9) at (4,0) [circle,draw]{$\ast$};
\draw (0)--(1)--(2)--(3);
\draw (4)--(5)--(6);
\draw (7)--(8)--(9);
\draw (1)--(4)--(7) (2)--(5)--(8) (3)--(6)--(9);
\draw[dashed] (2)--(4)--(8)--(6)--(2);
\end{tikzpicture}
\caption{The multiples of $8$}
\label{fig:8beishushu}
\end{subfigure}
\caption{The multiples of $6$ and $8$ on the nine-palace diagram}
\label{fig:68beishushu}
\end{figure}

Again, using the rotation invariance of the nine-palace diagram regarding multiples, the law of other even multiples is similar to the above situation: finding the multiples of an  even number is equivalent to counting numbers facing this even number. For example, finding the multiples of $4$ needs to face $4$, counting $1234$ and $6789$ in the order of head, left, right, and foot ($5$ is zero, and can be counted or not). Note that the head of $4$ is the primitive left midpoint, the left of facing $4$ is the primitive bottom midpoint, and so on. See Fig.~\ref{fig:4beishushu}.
Similarly, to find multiples of $6$ is to count towards $6$, and to find multiples of $8$ is to count towards $8$, see Fig.~\ref{fig:68beishushu}.

To sum up, we have

\begin{thm}[The counting method to find the multiples of a midpoint]\label{koujue:shushu2}
Using the nine-palace diagram to find a multiple of a midpoint is equivalent to counting $1234$ and $6789$ when standing at the center point and facing the given midpoint in the slopping square shape in order of 'front, left, right and back'.
\end{thm}

The theorems ~\ref{koujue:shushu1} and~\ref{koujue:shushu2} show that the units digit of a product has a good rule in the nine-palace diagram and thus easy to find.
So does the counting method determine the carry of a product? The answer is yes.

Recall the carry rule for addition of two digits: comparing the sum's ones digit to the addend, no carry $1$ is needed if the sum is larger (i.e., the point goes forward) and it requires carry $1$ if the sum is smaller (i.e., the point backward).
This carry rule may be generalized to the long sum of some numbers.

To find the sum of more numbers, you may add up the first two numbers, and then you start with the third number and then you add the subsequent numbers one by one. Observe the advance and retreat situation and carry $1$ for each backward time. Finally, look at the total number of backward times, and carry the number of backward times. For example. $3 + 6 + 7 + 8 + 2$ needs carry $2$,  because the total back times is $2 $: $3 + 6 = 9$ (forward), $9 + 7 = 16$ (backward), $6 + 8 = 14 $ (backward), $4 + 2 = 6 $ (forward).
In short, for the addition of a finite number of numbers, the times of going backwards is just the carry number. Note that the terms 'forward' or 'backward' are in the order of the size of the numbers in the primitive nine-palace diagram.

Since multiplication can be regarded as a special addition, it can be carried in the same way as addition, that is, the number of backward points is the number being carried. For example, $3\times 7$ carries $2$ because $3\times 7=7+7+7$ drops twice.

When we multiply a number by the counting method, it's easy to see how many times we must walk backwards. Therefore, we can use the counting method not only to compute the units digit of the product, but also to obtain the carry digit of this product. For example, in order to calculate $3\times 7$, looking at Fig.~\ref{fig:7beishushu} about the multiples of $7$, we count $1$ times, $2$ times and  $3$ times, and simultaneously observe the advance and retreat of the points. For example, in order to count $8\times 7$, again we look at the picture ~\ref{fig:7beishushu} and count numbers from $1$ to $8$ in turn and observe the advance and retreat of the points at the same time. It is easy to see that there are five retractions. To calculate $6\times 8$, look at the multiples of $8$~\ref{fig:8beishushu} and count from $1$ to $6$ in turn, while watching the advance and retreat of the points. It's easy to see that there are four retractions.
Summarizing the above methods, we get the following formula:

\begin{prop}[Find the carry digit of multiples by the counting method]\label{prop:shushufa-jinwei}
To find the product of two numbers using the counting method, the times of the points going backward  in the nine-palace diagram is just the carry digit.
\end{prop}

For quick calculations, it is essential to be familiar with the  carry situation of the product $9\times b$ on the graph, because it has more times of retreat than any product $a\times b$, and the carry situation of the latter can be seen on the carry situation graph of the former. Therefore, a carry diagram of $9\times b$ is actually the \emph{carry diagram} of the multiples of $b$.

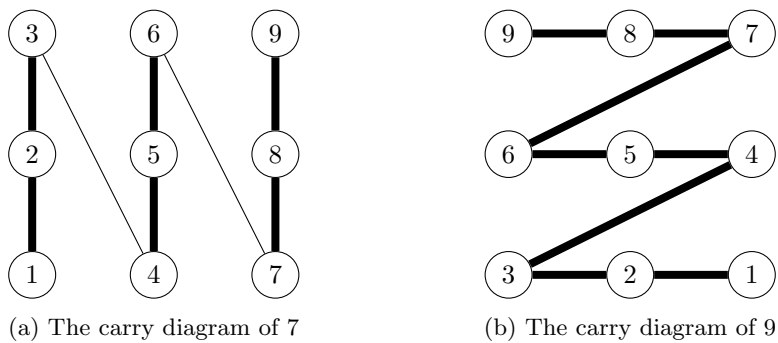
\begin{figure}[!htbp] 
\centering
\begin{subfigure}[t]{.4\linewidth}
\centering
\begin{tikzpicture}[scale=0.8]
\node(1) at (0,4) [circle,draw]{3};
\node(2) at (2,4) [circle,draw]{6};
\node(3) at (4,4) [circle,draw]{9};
\node(4) at (0,2) [circle,draw]{2};
\node(5) at (2,2) [circle,draw]{5};
\node(6) at (4,2) [circle,draw]{8};
\node(7) at (0,0) [circle,draw]{1};
\node(8) at (2,0) [circle,draw]{4};
\node(9) at (4,0) [circle,draw]{7};
\draw (1)--(8)(2)--(9);
\draw[line width=3pt] (1)--(4)--(7) (2)--(5)--(8) (3)--(6)--(9);
\end{tikzpicture}
\caption{The carry diagram of $7$}
\label{fig:7beishushu-jinwei}
\end{subfigure}
\quad
\begin{subfigure}[t]{.4\linewidth}
\centering
\begin{tikzpicture}[scale=0.8]
\node(1) at (0,4) [circle,draw]{9};
\node(2) at (2,4) [circle,draw]{8};
\node(3) at (4,4) [circle,draw]{7};
\node(4) at (0,2) [circle,draw]{6};
\node(5) at (2,2) [circle,draw]{5};
\node(6) at (4,2) [circle,draw]{4};
\node(7) at (0,0) [circle,draw]{3};
\node(8) at (2,0) [circle,draw]{2};
\node(9) at (4,0) [circle,draw]{1};
\draw[line width=3pt] (1)--(2)--(3);
\draw[line width=3pt] (4)--(5)--(6);
\draw[line width=3pt] (7)--(8)--(9);
\draw[line width=3pt] (3)--(4)(6)--(7);
\end{tikzpicture}
\caption{The carry diagram of $9$}\label{fig:9beishushu-jinwei}
\end{subfigure}
\caption{Carry diagrams of the multiples of $7$ and $9$ in the nine-palace diagram}
\label{fig:79beishushu-jinwei}
\end{figure}

Let's look at the carry diagram of a corner point $b$.
For example, the carry schematic diagram of multiples of $7$ and $9$ in the nine-palace diagram is shown in  Fig.~\ref{fig:79beishushu-jinwei}, in which the labels from $1$ to $9$ are successively connected by lines with thin lines indicating going forward and thick lines indicating going backward.
Fig.~\ref{fig:7beishushu-jinwei} has $6$ thick lines, i.e., total $6$ times of retreat, which means the carry of $9\times 7$ is $6$;  Fig.~\ref{fig:9beishushu-jinwei} has   $8$ thick lines, i.e.,  total $8$ times of retreat, which means that the carry of $9\times 9$ is $8$.
Similarly, you can sketch multiples of $1$ and multiples of $3$.

Next, let's look at the carry diagram of the midpoint $b$. The carry schematic diagram of multiples of $6$ and $8$ on the nine-palace diagram is shown in the Fig.~\ref{fig:68beishushu-jinwei}, where the dotted lines represent the primitive nine-palace diagram (which can be ignored) and the solid lines are connected from the label $1$ to $9$ in turn. The thin solid lines represent going forward and the thick lines represent going backward.
Fig.~\ref{fig:6beishushu-jinwei} has $5$ thick lines, i.e., $5$ times going backward, which means that $9\times 6$ needs carry $5$; In  Fig.~\ref{fig:8beishushu-jinwei}, there are $7$ thick lines, i.e., $7$ times going back, which means that the carry of  $9\times 8$ is $7$.
Similarly, you can plot multiples of $2$ and multiples of $4$ in the nine-palace diagram.

\begin{figure}[!htbp] 
\centering
\begin{subfigure}[t]{.4\linewidth}
\centering
\begin{tikzpicture}[scale=0.8]
\draw[step=2][dashed] (0,0) grid (4,4);
\node(0) at (-2,4) [circle,draw]{$5$};
\node(2) at (2,4) [circle,draw]{2,7};
\node(4) at (0,2) [circle,draw]{4,9};
\node(6) at (4,2) [circle,draw]{1,6};
\node(8) at (2,0) [circle,draw]{3,8};
\draw[line width=3pt] (4)..controls (1.3,1.3)..(8);
\draw[line width=3pt] (6)..controls (3.3,3.3)..(2);
\draw (2)..controls (2.3,2)..(8);
\draw (2)--(8) (0)--(6);
\draw[line width=3pt] (6)--(2)(8)--(4)--(0);
\end{tikzpicture}
\caption{The carry diagram of $6$}
\label{fig:6beishushu-jinwei}
\end{subfigure}
\quad
\begin{subfigure}[t]{.4\linewidth}
\centering
\begin{tikzpicture}[scale=0.8]
\draw[step=2][dashed] (0,0) grid (4,4);
\node(0) at (-2,4) [circle,draw]{$5$};
\node(2) at (2,4) [circle,draw]{4,9};
\node(4) at (0,2) [circle,draw]{3,8};
\node(6) at (4,2) [circle,draw]{2,7};
\node(8) at (2,0) [circle,draw]{1,6};
\draw[line width=3pt] (8)..controls (3.3,0.7)..(6);
\draw[line width=3pt] (6)..controls (2,1.7)..(4);
\draw[line width=3pt] (4)..controls (1.3,2.7)..(2);
\draw[line width=3pt](8)--(6)--(4)--(2)--(0);
\draw (0)..controls (-0.75,1.25)..(8);
\end{tikzpicture}
\caption{The carry diagram of $8$}
\label{fig:8beishushu-jinwei}
\end{subfigure}
\caption{Carry diagrams of the multiples of $6$ and $8$ on the nine-palace diagram}
\label{fig:68beishushu-jinwei}
\end{figure}

The carry schematic diagram for the center point $5$ consists of $5$ thick lines and $5$ thin lines connecting points $0$ and $5$, where the point $0$ corresponds    the  $0,2,4,6,8 $ times of $5$ and  the point $5$ corresponds  the  $1,3,5,7,9$ times of $5$. The reader can easily draw this graph by himself.

To calculate the carry digit of $a\times b$, we use the carry diagram for $b$. In fact, the graph consists mainly of a path from label $1$ to $9$. We start from the point marked $1$ and follow the path to the point marked $a$. We go through  total $a-1$ edges. By counting how many of the bold lines there are, we know how many retreats the points have and how many times of the carry the product $a\times b$ needs.
For example, the carry digit of $6\times 7$ is $4$, because the path from $1$ to $6$ in the carry diagram Fig.~\ref{fig:7beishushu-jinwei} of $9\times 7$ goes through  total $4$ bold edges.

Because of the particularity of the center point $5$ in the nine-palace diagram, finally we give the special method of the numeration method for its multiples. \par
According to Prop.~\ref{prop:5timesgewei}, the characteristic of the  units digit of $5$ times of each point on the primitive diagram is that the midpoints are attributed to $0$ and the corner points are attributed to  $5$. The corresponding carry rule is to carry one more digit for every even number, so the only points that affect the carry are the four midpoints, which are marked as thicker circles and exactly connected to a slopping square-shaped track, as shown in  Fig.~\ref{fig:5beishushu} with the dotted lines.

For example, to find the product $7\times 5$, we can find the point $7$ on the primitive nine-palace diagram. Because it is a corner point, the units digit of the product is equal to $5$. In the original nine-palace diagram, the number starts at $1$ and goes up to $7$ in sequence, so that there are three midpoints (i.e., the three larger circles in  Fig.~\ref{fig:7times5}), and thus the carry digit of the product is equal to $3$. Therefore, $7\times 5= 35$.
It can be seen that the multiples of $5$ can be counted in the original nine-palace diagram facing $1$ and use the diagram of the ones and carry rules of the multiples of $5$ on the nine-palace diagram ~\ref{fig:5beishushu}.

\begin{figure}[!htbp] 
\centering
\begin{subfigure}[t]{.4\linewidth}
\centering
\begin{tikzpicture}[scale=0.8]
\node(1) at (0,4) [circle,draw]{5};
\node(2) at (2,4) [circle,draw,line width=3pt]{0};
\node(3) at (4,4) [circle,draw]{5};
\node(4) at (0,2) [circle,draw,line width=3pt]{0};
\node(5) at (2,2) [circle,draw]{5};
\node(6) at (4,2) [circle,draw,line width=3pt]{0};
\node(7) at (0,0) [circle,draw]{5};
\node(8) at (2,0) [circle,draw,line width=3pt]{0};
\node(9) at (4,0) [circle,draw]{5};
\draw (1)--(2)--(3);
\draw (4)--(5)--(6);
\draw (7)--(8)--(9);
\draw (1)--(4)--(7) (2)--(5)--(8) (3)--(6)--(9);
\draw[dashed] (2)--(4)--(8)--(6)--(2);
\end{tikzpicture}
\caption{The schematic diagram of the units digits and the carry digits of the multiples of $5$ in the nine-palace diagram}
\label{fig:5beishushu}
\end{subfigure}
\quad
\begin{subfigure}[t]{.4\linewidth}
\centering
\begin{tikzpicture}[scale=0.8]
\node(1) at (0,4) [circle,draw]{1};
\node(2) at (2,4) [circle,draw,line width=3pt]{2};
\node(3) at (4,4) [circle,draw]{3};
\node(4) at (0,2) [circle,draw,line width=3pt]{4};
\node(5) at (2,2) [circle,draw]{5};
\node(6) at (4,2) [circle,draw,line width=3pt]{6};
\node(7) at (0,0) [circle,draw]{7};
\node(8) at (2,0) [circle,draw]{$\ast$};
\node(9) at (4,0) [circle,draw]{$\ast$};
\draw (1)--(2)--(3);
\draw (4)--(5)--(6);
\draw (7)--(8)--(9);
\draw (1)--(4)--(7) (2)--(5)--(8) (3)--(6)--(9);
\draw[dashed] (2)--(4)--(8)--(6)--(2);
\end{tikzpicture}
\caption{For example, $5\times 7$ by the nine-palace diagram}
\label{fig:7times5}
\end{subfigure}
\caption{The counting method for the multiple $5$ in the nine-palace diagram}
\label{fig:5beishushufa}
\end{figure}

\section{Simplification of the multiplicative carry diagram}

As we can see from the last section, although the multiplication carry is very regular, we need to count the number of backward edges, which is not convenient in practical application. For this reason, we shift the perspective from the edge to the point, so that the points that need to carry are drawn as black, and the points that don't need to carry are drawn as white. As a consequence,  the carry diagram can be greatly simplified, and it is very convenient to use.

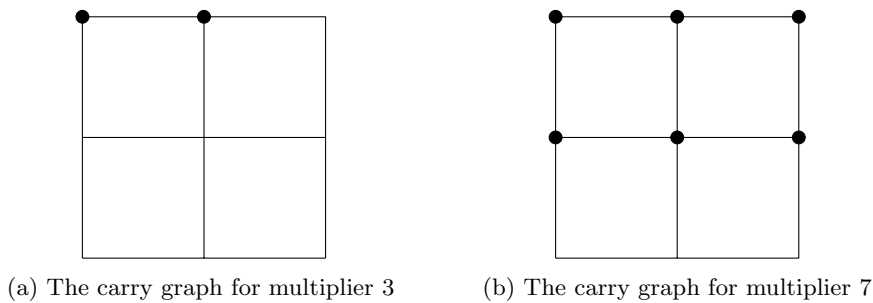
\begin{figure}[!htbp] 
\centering
\begin{subfigure}[t]{.4\linewidth}
\centering
\begin{tikzpicture}[scale=0.8]
\draw[step=2] (0,0) grid (4,4);
\draw[fill] (0,4) circle (3pt)(2,4) circle (3pt);
\end{tikzpicture}
\caption{The carry graph for multiplier $3$}
\label{fig:3-7jinweiJH-a}
\end{subfigure}
\quad
\begin{subfigure}[t]{.4\linewidth}
\centering
\begin{tikzpicture}[scale=0.8]
\draw[step=2] (0,0) grid (4,4);
\foreach \x in {0,2,4} \foreach \y in {2,4} \draw[fill] (\x,\y) circle (3pt);
\end{tikzpicture}
\caption{The carry graph for multiplier $7$}
\label{fig:3-7jinweiJH-b}
\end{subfigure}
\caption{The carry graph for multipliers $3$ and $7$}
\label{fig:3-7jinweiJH}
\end{figure}

The simplified carry diagram of the multiplier $3$ and $7$ is shown in ~\ref{fig:3-7jinweiJH}.  By comparing the Fig.~\ref{fig:3-7jinweiJH-b} with  Fig.~\ref{fig:9beishushu-jinwei}, it can be seen that the former is much simpler than the latter, although both are carry diagrams for multiples of $7$.
In  Fig.~\ref{fig:3-7jinweiJH-b}, when counting numbers facing the primitive point $7$, no new carry generates when counting at $1$, $4$  and $7$, so these three points are not marked black, while the other six points will generate new carry, so they are all marked with black dots.

When the counting reaches $2$, the carry is $1$, because  $2$  is the first black dot we encountered. When $3$ is reached, $2$ is carried, because the second black dot is encountered. When the counting reaches $5$, the carry is $3$, because the third black dot is encountered. When $6$ is reached, $4$ is carried, because the fourth black dot is encountered. When the counting reaches $8$, the carry is $5$, because the fifth black dot is encountered. When the counting reaches $9$, the carry is $6$, because the sixth black dot is encountered. When $1$ is counted, there is no carry; When the counting reaches $4$, the rounding number does not increase and $2$ is still carried, because two black dots have already passed by this time. When the counting reaches $7$, the rounding number does not increase, but $4$ is still carried, because four black dots have already been passed.  We see that the graph ~\ref{fig:3-7jinweiJH-b} has carry dots of the total number $6$. One can understand similarly the carry diagram~\ref{fig:3-7jinweiJH-b} for the multiplier $3$, where there are exactly $2$ carry points.

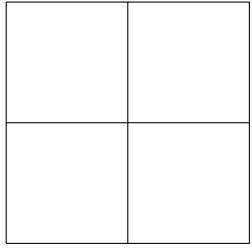
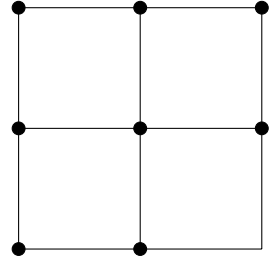
\begin{figure}[!htbp] %
\centering
\begin{subfigure}[t]{.4\linewidth}
\centering
\begin{tikzpicture}[scale=0.8]
\draw[step=2] (0,0) grid (4,4);
\end{tikzpicture}
\caption{The carry diagram of the multiplier $1$ }
\label{fig:1-9jinweiJH-a}
\end{subfigure}
\quad
\begin{subfigure}[t]{.4\linewidth}
\centering
\begin{tikzpicture}[scale=0.8]
\draw[step=2] (0,0) grid (4,4);
\foreach \x in {0,2,4} \foreach \y in {2,4} \draw[fill] (\x,\y) circle (3pt);
\draw[fill] (0,0) circle (3pt)(2,0) circle (3pt);
\end{tikzpicture}
\caption{The carry diagram of the multiplier $9$}
\label{fig:1-9jinweiJH-b}
\end{subfigure}
\caption{The carry diagram of the multipliers $1$ and $9$}
\label{fig:1-9jinweiJH}
\end{figure}

The carry diagram of the multiplier $1$ and $9$ is shown in Fig.~\ref{fig:1-9jinweiJH}, where there is no carry point for multiplier $1$, while the other eight points except the point 9 are the carry points for multiplier 9.

\begin{figure}[!htbp] 
\centering
\begin{tikzpicture}[scale=0.8]
\draw[step=2] (0,0) grid (4,4);
\draw[fill] (0,2) circle(3pt)(2,0) circle(3pt)(2,4) circle(3pt)(4,2) circle(3pt);
\end{tikzpicture}
\caption{The carry diagram of multiplier $5$}
\label{fig:5jinweiJH}
\end{figure}
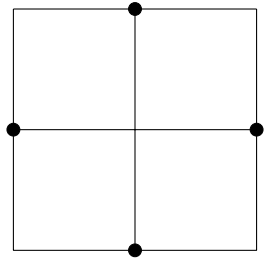

The carry diagram of the multiplier $5$ is slight different from that described in the previous section. See Fig.~\ref{fig:5jinweiJH}. You need to count numbers facing $1$. Obviously, each of $2$, $4$, $6$, and $8$ adds a carry, so the midpoints are all black dots. We know that the units digit of an even number multiplied by $5$ is $0$, and the units digit of an odd number multiplied by $5$ is $5$. Thus, if the multiplicative units digit is obtained from this carry map, the midpoint and zero become to the zero $0$, and the corner points and  the center point  become  to the  center point $5$.

\begin{figure}[!htbp] 
\centering
\begin{subfigure}[t]{.4\linewidth}
\centering
\begin{tikzpicture}[scale=0.8]
\draw[step=2] (0,0) grid (4,4);
\draw[fill] (-2,4) circle(3pt);
\draw (-2,4)--(0,4);
\end{tikzpicture}
\caption{The carry diagram of multiplier $2$}
\label{fig:2-4jinweiJH-a}
\end{subfigure}
\quad
\begin{subfigure}[t]{.4\linewidth}
\centering
\begin{tikzpicture}[scale=0.8]
\draw[step=2] (0,0) grid (4,4);
\draw[fill] (-2,4) circle(3pt)(2,4) circle(3pt);
\draw (-2,4)--(0,4);
\end{tikzpicture}
\caption{The carry diagram of multiplier $4$}
\label{fig:2-4jinweiJH-b}
\end{subfigure}
\caption{The carry diagram of multipliers $2$ and $4$}
\label{fig:2-4jinweiJH}
\end{figure}
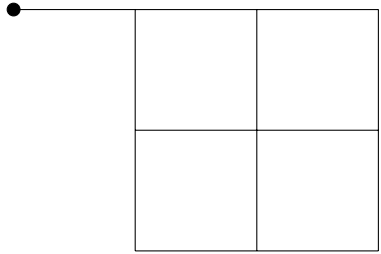
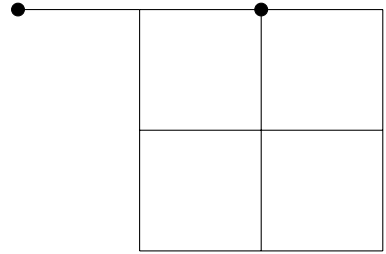

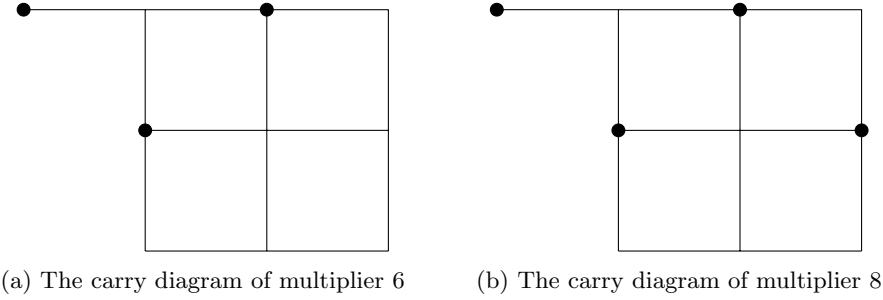
\begin{figure}[!htbp] 
\centering
\begin{subfigure}[t]{.4\linewidth}
\centering
\begin{tikzpicture}[scale=0.8]
\draw[step=2] (0,0) grid (4,4);
\draw[fill] (-2,4) circle(3pt)(0,2) circle(3pt)(2,4) circle(3pt);
\draw (-2,4)--(0,4);
\end{tikzpicture}
\caption{The carry diagram of multiplier $6$}
\label{fig:6-8jinweiJH-a}
\end{subfigure}
\quad
\begin{subfigure}[t]{.4\linewidth}
\centering
\begin{tikzpicture}[scale=0.8]
\draw[step=2] (0,0) grid (4,4);
\draw[fill] (-2,4) circle(3pt)(0,2) circle(3pt)(2,4) circle(3pt)(4,2) circle(3pt);
\draw (-2,4)--(0,4);
\end{tikzpicture}
\caption{The carry diagram of multiplier $8$}
\label{fig:6-8jinweiJH-b}
\end{subfigure}
\caption{The carry diagram of multipliers $6$ and $8$}
\label{fig:6-8jinweiJH}
\end{figure}

When the multiple of $n$ is even (the midpoint), the total number of carry points is also equal to $n-1$. However, we can only see just half of the black dots in the figure. This is because except for zero, the other black dots need to be counted twice.
See Fig~\ref{fig:2-4jinweiJH} and Fig.~\ref{fig:6-8jinweiJH}. We can compare the carry diagram ~\ref{fig:6-8jinweiJH-a} of the multiples of $6$  with  Fig.~\ref{fig:6beishushu-jinwei} in the previous section.
Counting numbers towards the primitive point $6$, there is no carry when counting $1$; When the counting numbers reaches $2$, we have carry $1$, because that is the first black dot encountered. When $3$ is counted, no new carry is generated, and now $1$ is still carried, because a black dot has passed by this time. When $4$ is reached, we have carry $2$, because the second black dot is encountered. When the counting  reaches $5$, we need carry $3$, because the third black dot is encountered. When $6$ is counted, no new carry is generated, and now $3$ is still carried, because three black dots have passed by this time. When $7$ is reached, we carry $4$, because the fourth black dot is encountered (one of the black dots must be counted twice); When $8$ is counted, no new carry is generated, and $4$ is still carried, because four black dots have been passed by this time. When $9$ is reached, carry $5$, because the fifth black dot is encountered (two of which must be counted twice).

From the  above discussion,  we obtain the following conclusion:

\begin{thm}[The carry theorem in the nine-palace diagram]
There are $n-1$ carry points, which are exactly the points from $1$ to $9$ that are less than the primitive point $n$.
\end{thm}

This theorem tells us that the carry of multiplication also shows good regularity in the nine-palace diagram, so it is also easy to grasp.

\section{The counting method in the nine-palace diagram are used for multiplication and division}

Since the units and carry rules of multiplication are well represented in the nine-palace diagram, which can be used to fast multiplication and division. The multiplication and division of a long number by a one-digit number, except for carry and abdication, are simply a rotation in the nine-palace diagram. If you can multiply and divide a long number by a one-digit number, you can multiply and divide  a long number by a long number.

Let the multiplicand $a_1a_2\ldots a_n$ be a $n$-bit number and the multiplier $b$ be a one-bit number. To compute the product $a_1a_2\ldots a_n\times b$, we consider the product $a_i\times b$ in turn, where $i =1,2,\ldots,n$.
Set $a_i\times b=(J_i,g_i)$. Then we are dealing with the  $i$th bit (called the current bit), which units digit $g_i$ is called the \emph{current units digit} and which carry digit $J_i$ is called the \emph{current carry}. The current carry $J_{i+1}$ of the next bit is called the \emph{current back carry} or briefly, the \emph{back carry} .  Define
\begin{equation}
\textrm{the {\it bit product}$:=$the current units digit $+$ the back carry}.
\end{equation}
For $i=0,1,2,\ldots,n$, denoting the $i$th bit product as $p_i$, we have
\begin{equation}\label{eq:benge-jinwei}
p_i=g_i+J_{i+1},
\end{equation}
Here we convent $g_0=J_{n+1}=0$. Then we obtain
\begin{displaymath}
\begin{array}{rl}
a_1a_2\ldots a_n\times b &=p_0p_1\ldots p_n\\
                         &=g_0g_1\ldots g_n+J_1J_2\ldots J_{n+1}\\
                         &=0g_1\ldots g_n+J_1J_2\ldots J_n0.
\end{array}
\end{displaymath}
We call $0g_1\ldots g_n$ the \emph{units digit sequence}, and call $J_1J_2\ldots J_n0$  the \emph{carry sequence}. The above conclusion can be written as the following
\begin{prop}
When the multiplier is a single digit, the product is the sum of the units digit sequence and the carry sequence.
\end{prop}

Both the units digit sequence and the carry sequence can be easily seen from the nine-palace diagram using the counting method. Therefore, as long as we can add two long numbers  on the nine-palace diagram, we can complete the multiplication of a long number and a one-digit number. The specific method is summarized as the following

\begin{prop}[The multiplication with the multiplier a units digit]
When the multiplier is a units digit, the direction is determined by the multiplier and the units digit sequence and the carry sequence are read out according to the counting method. Then the sum of the two sequences is just the  required product, which can be read out directly in the primitive nine-palace diagram.
\end{prop}

\begin{exmp}\rm
Use the counting method to calculate the produc $4\,789\times 3.$
\end{exmp}

\begin{solution}
Since the multiplier is $3$, we face the right side of the nine-palace diagram  (where the position $3$ is pronounced as $1$) and read out the multiplicand $04\,789$, which is denoted by $JABCD$. According to the formula~\ref{koujue:shushu1}, this is actually the units digit sequence, as shown in  Fig.~\ref{fig:4789*3-a}.

\begin{figure}[!htbp] 
\centering
\begin{subfigure}[t]{.4\linewidth}
\centering
\begin{tikzpicture}[scale=0.6]
\draw[step=2] (0,0) grid (4,4);
\draw[fill] (0,0) circle (3pt)  (0,2) circle (3pt)  (0,4) circle (3pt) (2,4) circle (3pt)(4,6) circle (3pt);
\coordinate [label=180:$J$](J) at (4,6);
\coordinate [label=135:$A$](A) at (2,4);
\coordinate [label=135:$B$](B) at (0,4);
\coordinate [label=180:$C$](C) at (0,2);
\coordinate [label=180:$D$](D) at (0,0);
\draw[line width=3pt,->,>=latex](J)--(A)--(B)--(C)--(D);
\end{tikzpicture}
\caption{The units digit sequence}
\label{fig:4789*3-a}
\end{subfigure}
\quad
\begin{subfigure}[t]{.4\linewidth}
\centering
\begin{tikzpicture}[scale=0.6]
\draw[step=2] (0,0) grid (4,4);
\draw[fill] (0,4) circle (3pt)(4,6) circle (3pt);
\draw[fill] (2,4) circle (3pt)[shift={(-4pt,4pt)}](2,4) circle (3pt);
\draw[fill] [shift={(4pt,-4pt)}](2,4) circle (3pt);
\coordinate [label=180:$J$](J) at (0,4);
\coordinate [label=135:$A$](A) at (2,4);
\coordinate [label=-45:$B$](B) at (2,4);
\coordinate [label=-135:$C$](C) at (2,4);
\coordinate [label=180:$D$](D) at (4,6);
\draw[line width=3pt][->,>=latex](J)--(B)--(D);
\end{tikzpicture}
\caption{The carry consequence}
\label{fig:4789*3-b}
\end{subfigure}
\quad
\begin{subfigure}[t]{.4\linewidth}
\centering
\begin{tikzpicture}[scale=0.6]
\draw[step=2] (0,0) grid (4,4);
\draw[fill] (0,0) circle (3pt)  (0,2) circle (3pt)  (0,4) circle (3pt) (4,4) circle (3pt)(4,2) circle (3pt);
\coordinate [label=180:$J$](J) at (0,4);
\coordinate [label=135:$A$](A) at (0,2);
\coordinate [label=135:$B$](B) at (4,4);
\coordinate [label=135:$C$](C) at (4,2);
\coordinate [label=180:$D$](D) at (0,0);
\draw[line width=3pt,->,>=latex](J)--(A)--(B)--(C)--(D);
\end{tikzpicture}
\caption{The production}
\label{fig:4789*3-c}
\end{subfigure}
\caption{Use the counting method to calculate the product $4\,789\times 3$}
\label{fig:4789*3}
\end{figure}
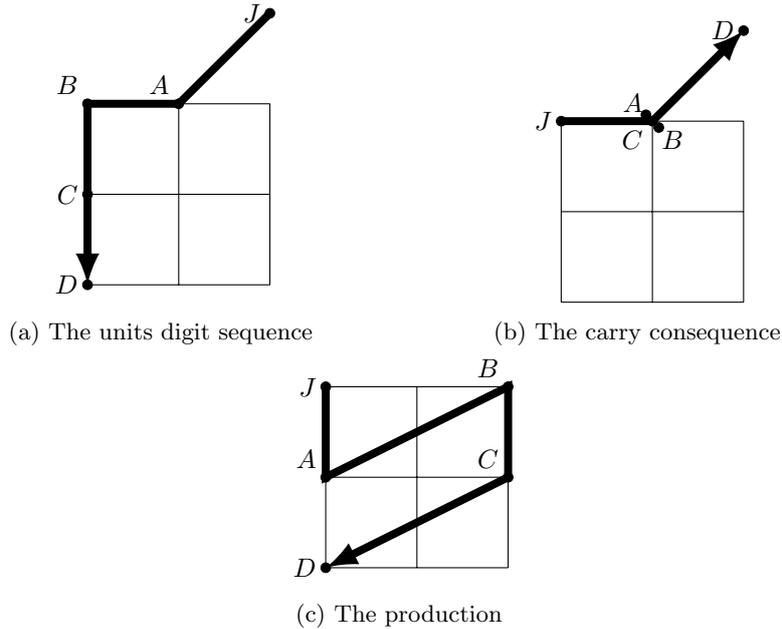
\par
According to the formula~\ref{prop:shushufa-jinwei}, one can count out the carry consequence in Fig.~\ref{fig:4789*3-a}. For example, the points $A,B$ correspond to carry digits $1,2$, respectively, because there are 1 and 2 times of going backward from 1 to points $A,B$. Thus we get the carry consequence, denoted still as $JABCD$, see Fig.~\ref{fig:4789*3-b}.
\par
Sum up the points denoted as the same letters in Figs.~\ref{fig:4789*3-a} and \ref{fig:4789*3-b}, so as to obtain the product consequence, denoted still as $JABCD$, see Fig.~\ref{fig:4789*3-c}.
In the primitive nine-palace diagram, points $JABCD$ are read as $14367$. Thus $4\,789\times 3=14\,367$.
\end{solution}

\begin{exmp}\rm
Use the counting method to calculate the produc $92\,867\times 8.$
\end{exmp}

\begin{solution}
Since the multiplier is $3$, we face the bottom side of the nine-palace diagram  (where the position $9$ is pronounced as $1$) and read out the multiplicand $092\,867$, which is denoted by $J$, $\ddot{A}$, $\dot{B}$, $\ddot{C}$, $\ddot{D}$, $\ddot{E}$, where the superscript dots indicate that the corresponding points are located in the first or second circle of the flopping square shape formed by the four middle points when counting numbers on this square. According to Prop.~\ref{koujue:shushu1},  this is actually the units digit sequence, as shown in  Fig.~\ref{fig:92867*8-a}.

\begin{figure}[!htbp] 
\centering
\begin{subfigure}[t]{.4\linewidth}
\centering
\begin{tikzpicture}[scale=0.6]
\draw[step=2] (0,0) grid (4,4);
\draw[fill] (0,2) circle (3pt)  (4,6) circle (3pt)  (2,0) circle (3pt)  (2,4) circle (3pt) (4,2) circle (3pt)[shift={(4pt,4pt)}](4,2) circle (3pt);
\coordinate [label=180:$J$](J) at (4,6);
\coordinate [label=135:$\ddot{A}$](A) at (2,4);
\coordinate [label=45:$\dot{B}$](B) at (4,2);
\coordinate [label=180:$\ddot{C}$](C) at (0,2);
\coordinate [label=-90:$\ddot{D}$](D) at (2,0);
\coordinate [label=-45:$\ddot{E}$](E) at (4,2);
\draw[line width=3pt,->,>=latex](J)--(A)--(B)--(C)--(D)--(E);
\end{tikzpicture}
\caption{The units consequence}
\label{fig:92867*8-a}
\end{subfigure}
\quad
\begin{subfigure}[t]{.4\linewidth}
\centering
\begin{tikzpicture}[scale=0.6]
\draw[step=2] (0,0) grid (4,4);
\draw[fill] (0,0) circle (3pt)(0,2) circle (3pt)(0,4) circle (3pt)(2,2) circle (3pt)(4,2) circle (3pt)(4,6) circle (3pt);
\coordinate [label=180:$J$](J) at (0,0);
\coordinate [label=180:$A$](A) at (0,4);
\coordinate [label=0:$B$](B) at (4,2);
\coordinate [label=180:$C$](C) at (0,2);
\coordinate [label=-135:$D$](D) at (2,2);
\coordinate [label=0:$E$](E) at (4,6);
\draw[line width=3pt][->,>=latex](J)--(A)--(B);
\draw[line width=3pt][->,>=latex](C)--(D)--(E);
\end{tikzpicture}
\caption{The carry consequence}
\label{fig:92867*8-b}
\end{subfigure}
\quad
\begin{subfigure}[t]{.4\linewidth}
\centering
\begin{tikzpicture}[scale=0.6]
\draw[step=2] (0,0) grid (4,4);
\draw[fill] (0,0) circle (3pt)  (0,2) circle (3pt)  (2,4) circle (3pt) (4,0) circle (3pt)(4,2) circle (3pt)(4,4) circle (3pt);
\coordinate [label=180:$J$](J) at (0,0);
\coordinate [label=135:$A$](A) at (0,2);
\coordinate [label=90:$B$](B) at (2,4);
\coordinate [label=0:$C$](C) at (4,0);
\coordinate [label=0:$D$](D) at (4,4);
\coordinate [label=0:$E$](E) at (4,2);
\draw[line width=3pt,->,>=latex](J)--(A)--(B)--(C);
\draw[line width=3pt,->,>=latex](D)--(E);
\end{tikzpicture}
\caption{The production}
\label{fig:92867*8-c}
\end{subfigure}
\caption{Using the counting method to calculate the product $92\,867\times 8$}
\label{fig:92867*8}
\end{figure}
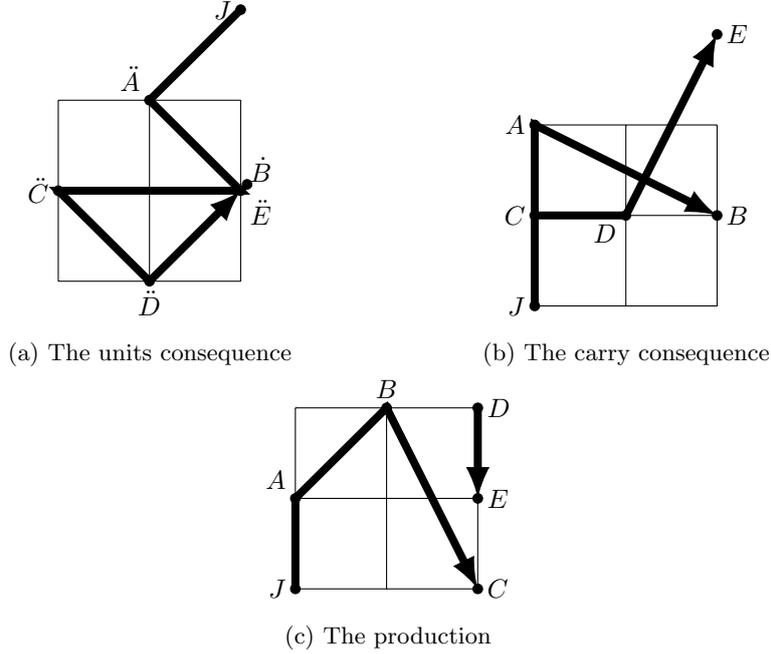
\par
According to Prop.~\ref{prop:shushufa-jinwei}, one can read out the carry consequence from  Fig.~\ref{fig:92867*8-a}. For example, the points $\ddot{A},\dot{B}$ correspond to the carry digits $7,1$, respectively, since there will be 7 and 1 times  of abdication when facing down and counting from 1 to points $A,B$. It is noted that $\ddot{A}$ is located at the second round of the flopping square and that $\dot{B}$ is located at the first round of the square.  $A$ has carry digit $J$, $B$ has carry digit $A$, $C$ has carry digit $B$, and so on. Thus we obtain the carry consequence $JABCDE$, as shown in Fig.~\ref{fig:92867*8-b}.
\par
Sum up the points denoted as the same letters in Figs.~\ref{fig:92867*8-a} and~\ref{fig:92867*8-b},  so as to obtain the product consequence, denoted still as $JABCDE$, as shown in Fig.~\ref{fig:92867*8-c}.
In the primitive nine-palace diagram, points $JABCD$ are read as $742936$. Thus $92\,867\times 8=742\,936$.
\end{solution}

Transposing the terms of the formula~\ref{eq:benge-jinwei}, we obtain immediately: for any $i=0,1,2,\ldots,n$,
\begin{equation}\label{eq:benge-jinwei-1}
g_i=p_i-J_{i+1}.
\end{equation}
Since this formula can be used to division, we put it in the following proposition.

\begin{prop}
The current units digit equals to the difference of the bit product and the back carry (sometimes probably minus 1).
\end{prop}

Next,  we show how to use this formula for a division with a single digit divisor (in the exact divisibility case).
First, since $J_{n+1}=0$,  we have $g_n=p_n$. Second, we can start from point $p_n$ and read out $J_n$ by counting method on the nine-palace diagram. Thus, $g_{n-1}=p_{n-1}-J_n$. Note that if the difference is negative, we must borrow 1 from $p_{n-2}$ as 10 in the current bit, that is, we obtain $p_{n-2}-1$ and $p_{n-1}+10$. Use the new value of $p_{n-1}$ to make difference $g_{n-1}=p_{n-1}-J_n$, so as to obtain a non-negative number.
Continuing to consider the higher bit,  we have $g_{n-2}=p_{n-2}-J_{n-1}$, where $p_{n-2}$ is the newest, that is, it is possibly equals to the old value minus $1$. As above, it is noted that if the difference is negative, we must borrow 1 from the higher bit. Repeating this process, we obtain at last that $g_0=p_0-J_1=0$. Consequently, we obtain the units consequence $0g_1g_2\ldots g_n$, from which we can read out the quotient of the division by the counting method on the nine-palace diagram. So we have:

\begin{prop}[The division with the exact divisibility case]
In the  exact divisibility  case,  the division with the divisor a units digit can be made on the nine-palace diagram in the following way: if the divisor is the corner point, one first read the dividend in the primitive nine-palace diagram, then from the higher bit to the lower bit make the difference of the bit product and the back carry digit (maybe being minus 1 is needed) so as to obtain the current units digit of the quotient, from which one can read the units digit consequence by the counting method, which is the desired quotient.
\end{prop}

\begin{exmp}\rm
Use the counting method to make the division $14\,367\div 3.$
\end{exmp}

\begin{solution}
The divisor is $3$, or equivalently, the multiplier is $3$. So when using the counting method, we must face the right side (i.e., read the primitive position $3$ as $1$).
\par
In the primitive nine-palace diagram, we label the dividend consequence $14367$ as $JABCD$, which are connected by some lines as in Fig.~\ref{fig:4789div3-a}.
The last point $D$ stays still. According to the formula~\ref{prop:shushufa-jinwei}, we get that the point $D$ corresponds to the carry  $J_D=2$, since facing the right and counting numbers from $1$ to the point $D$ there is $2$ abdications, also see Fig.~\ref{fig:3-7jinweiJH-a},  the carry diagram of multiplier $3$. Make difference $C-J_D$, and resulting point is still labeled as $C$. For this new point $C$, by the counting method, corresponds the carry digit  $J_C=2$.
Make difference $B-J_C$, and the resulting point is still labeled as $B$. For this new point $B$, by the counting method again, it corresponds to the carry digit $J_B=2$.  Again, make the difference $A-J_B$, which is still labeled as $A$. For  this new point $A$, by the counting method again, it corresponds to the carry digit $J_A=1$. Again, make the difference $J-J_A$, which is still labeled as $J$. Since $J$ is at $0$, the operation is over. Now, we obtain the units digit consequence $JABCD$ shown as in Fig.~\ref{fig:4789div3-b}.
\begin{figure}[!htbp] 
\begin{subfigure}[t]{.4\linewidth}
\centering
\begin{tikzpicture}[scale=0.6]
\draw[step=2] (0,0) grid (4,4);
\draw[fill] (0,0) circle (3pt)  (0,2) circle (3pt)  (0,4) circle (3pt) (4,4) circle (3pt)(4,2) circle (3pt);
\coordinate [label=180:$J$](J) at (0,4);
\coordinate [label=135:$A$](A) at (0,2);
\coordinate [label=135:$B$](B) at (4,4);
\coordinate [label=135:$C$](C) at (4,2);
\coordinate [label=180:$D$](D) at (0,0);
\draw[line width=3pt,->,>=latex](J)--(A)--(B)--(C)--(D);
\end{tikzpicture}
\caption{The dividend}
\label{fig:4789div3-a}
\end{subfigure}
\quad
\centering
\begin{subfigure}[t]{.4\linewidth}
\centering
\begin{tikzpicture}[scale=0.6]
\draw[step=2] (0,0) grid (4,4);
\draw[fill] (0,0) circle (3pt)  (0,2) circle (3pt)  (0,4) circle (3pt) (2,4) circle (3pt)(4,6) circle (3pt);
\coordinate [label=180:$J$](J) at (4,6);
\coordinate [label=135:$A$](A) at (2,4);
\coordinate [label=135:$B$](B) at (0,4);
\coordinate [label=180:$C$](C) at (0,2);
\coordinate [label=180:$D$](D) at (0,0);
\draw[line width=3pt,->,>=latex](J)--(A)--(B)--(C)--(D);
\end{tikzpicture}
\caption{The units digit consequence}\label{fig:4789div3-b}
\end{subfigure}
\caption{Using the counting method to make the division $14\,367\div 3$}
\label{fig:4789div3}
\end{figure}
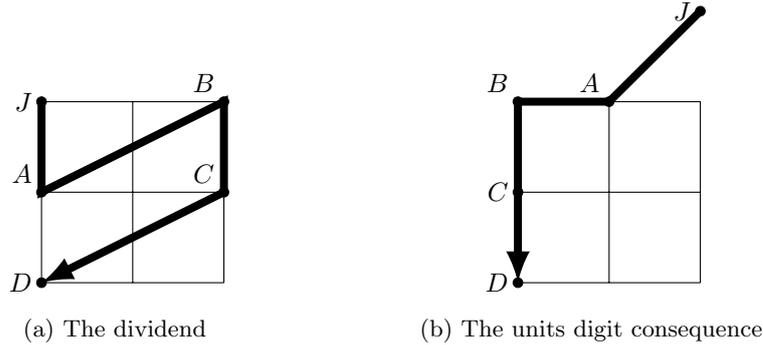
\par
We face the right of the nine-palace diagram and from  Fig.~\ref{fig:4789div3-b} read out the units digit consequence $JABCD$, so that we obtain $04\,789$, which means that $14\,367\div 3=4\,789.$
\end{solution}

The essence of the above method is to find the units digit sequence of the product of the quotient and the divisor. The characteristic of this method is to start the calculation from the lowest bit.  The method is elegant but has a lot of limitations. First, when the divisor is an even number of the one-digit, it does not work in this way when dividing on a nine-palace diagram, because different multiples may correspond to the same point when counting numbers toward the midpoint. Second, in the case of not being divisible, one neither can use the above method. Below we introduce a method of counting numbers from the top bit, which is generally applicable to one-digit division, regardless the divisor odd or even, and regardless whether or not it is exactly divisible:

\begin{prop}[The division with a single digit divisor]
For division with a single digit as the divisor,  the divisor-oriented counting method can be used, and the specific operation process is as follows:
\par
(1) Mark the dividend in the primitive nine-palace diagram. If the first digit of the dividend is not less than the divisor, add a zero in front of the dividend as the highest digit. \par
(2) Pay attention to the highest two digits  $XY$ of the dividend and determine a quotient point by counting. When the sequential counting reaches the quotient point, the total number of steps backward should be exactly equal to $X$, and the position of the quotient point (in the primitive nine-palace diagram) should be as far back as possible but not more than $Y$. Use the quotient point to subtract $Y$ to get the temporary remainder. \par
(3)  If the above requirement cannot be met, then reduce the number of backward steps by $1$. At this time, $10$ should be added after subtracting $Y$ with the quotient point, thus obtaining the temporary remainder. \par
(4) The temporarily remainder is the highest digit of the new dividend, the other digits of the new dividend is those that have not been involved in the previous calculation. Use the  new dividend instead of the old dividend to repeat the above processes (1)-(3) until the last units digit of the dividend, and then  the temporary remainder in the final step is just the remainder of the division; \par
(5) Read the sequence of quotient points as a divisor-oriented counting method, which is the desired quotient.
\end{prop}

\section{Conclusion}

As we have seen above, arithmetic operations show many beautiful patterns on the nine-palace diagram. Because the graph is intuitive and easy to remember, the nine-palace diagram is a good tool for mental arithmetic. Note that the  arithmetic by the nine-palace diagram is not the same as the usual finger-pointing algorithms or the China's traditional quick calculation method of wallowing gold up in your sleeve. The theory we have established in this paper can be called the nine-palace arithmetic, which is actually the geometric theory of arithmetic. Because this theory is actually based on Chinese traditional culture,  it has a certain philosophical significance.

The theory of the nine-palace arithmetic has a very rich contents. In addition to some aspects introduced in this paper, there are also some other contents such as the partition carry method, the straw man multiplication formulae, the multiplication track, the addition track, the addition road and so on. But limited by space, we do not make further discussion here.


\end{sloppypar}  

\end{document}